\documentclass[journal,twoside,web]{ieeecolor}
\usepackage{etoolbox}
\newtoggle{onecolumnversion}

\makeatletter
\@ifclasswith{ieeecolor}{onecolumn}{\toggletrue{onecolumnversion}}{\togglefalse{onecolumnversion}}
\makeatother

\usepackage[usenames,dvipsnames]{xcolor}
\usepackage{generic}
\usepackage[sort,noadjust]{cite}
\usepackage{amsmath,amssymb,amsfonts}
\usepackage{nicefrac}
\usepackage{caption} 
\iftoggle{onecolumnversion}{}{
\captionsetup{skip=0pt}
\setlength{\floatsep}{1.0pt plus 1.0pt minus 1.0pt}
\setlength{\textfloatsep}{1.0pt plus 1.0pt minus 1.0pt}
}

\usepackage{soul}
\sethlcolor{cyan}
\usepackage{textcomp}
\def\BibTeX{{\rm B\kern-.05em{\sc i\kern-.025em b}\kern-.08em
    T\kern-.1667em\lower.7ex\hbox{E}\kern-.125emX}}
\usepackage[utf8]{inputenc}
\usepackage[T1]{fontenc}
\usepackage[colorlinks=true,    %no frame around URL
      urlcolor=black,     %no colors
      citecolor=blue,
      menucolor=black,    %no colors
      linkcolor=black,    %no colors 
      breaklinks=true,
]{hyperref}
\usepackage{graphicx}

\usepackage{amssymb,amsmath}

\usepackage{xparse}
\usepackage{amsmath}
\usepackage{bm}
\usepackage{mathtools}
\mathtoolsset{centercolon}
 
\usepackage{tabularx}

\usepackage{cases}
\usepackage{diagbox}
\usepackage{algorithm}
\usepackage{algcompatible}
\usepackage{booktabs}
\usepackage{multirow}

\usepackage{lipsum}
\usepackage{tabularx}
\usepackage{tcolorbox}

\newcommand\numberthis{\addtocounter{equation}{1}\tag{\theequation}}

\newcommand{\Rl}[2]{\ensuremath{\mathbb{R}^{#1}_{#2}}}   % Euclidean space
\newcommand{\R}{\ensuremath{\mathbb{R}}}
\newcommand{\Rlp}{\ensuremath{\mathbb{R}_{>0}}}

\newcommand{\Rlo}{\ensuremath{\mathbb{R}_{\geq 0}}}

\newcommand{\Zo}{\ensuremath{\mathbb{Z}_{\geq 0}}}
\newcommand{\Zp}{\ensuremath{\mathbb{Z}_{> 0}}}
\definecolor{bleucit}{rgb}{0.2,0.4,0.6} % nouvelle couleur

\newcommand{\dst}{\displaystyle}

\definecolor{blue_cv}{rgb}{0.09,0.35,0.78}
\newcommand{\alphav}{\ensuremath{\overline{\alpha}_V}}

\newcommand{\KL}{\ensuremath{\mathcal{KL}}}
\newcommand{\K}{\ensuremath{\mathcal{K}}}
\newcommand{\Kinf}{\ensuremath{\mathcal{K}_{\infty}}}

\newcommand{\Argmin}{\ensuremath{\text{argmin}\,}}
\DeclareMathOperator*{\argmin}{\Argmin}

\newcommand{\dom}{\ensuremath{\text{dom}\,}}

\newtheorem{defn}{Definition}

\newtheorem{ass}{Assumption}

\newtheorem{prop}{Proposition}

\newtheorem{lem}{Lemma}

\newtheorem{thm}{Theorem}
\newtheorem{cor}{Corollary}

\newtheorem{rem}{Remark}
\newtheorem{sass}{Standing Assumption}

\usepackage[all]{hypcap}

\title{\LARGE \bf
Policy iteration: for want of recursive feasibility, all is not lost}

\author{
Mathieu Granzotto, Olivier  Lindamulage De Silva, Romain Postoyan, Dragan Nešić, and {Zhong-Ping Jiang}
\thanks{Mathieu Granzotto and Dragan Nešić are with the Department of Electrical and Electronic Engineering, University of Melbourne, Parkville, VIC 3010, Australia (e-mail: \{mgranzotto, dnesic\}@unimelb.edu.au). Their work was supported by the Australian Research Council under the Discovery Project DP210102600.}
\thanks{Olivier  Lindamulage De Silva and  Romain Postoyan  are with the Universit\'e de Lorraine, CNRS, CRAN, F-54000 Nancy, France (e-mails: \{name.surname\}@univ-lorraine.fr).}
\thanks{Zhong-Ping Jiang is with the Department of Electrical and Computer Engineering, New York University, 370 Jay Street, Brooklyn, NY 11201, USA. email: zjiang@nyu.edu. His work was supported partly by the National Science Foundation under Grant EPCN-1903781.}
\thanks{This work was supported by the France-Australia collaboration project IRP-ARS CNRS.}} %  stops a space

\makeatletter
\patchcmd{\hyper@makecurrent}{%
    \ifx\Hy@param\Hy@chapterstring
        \let\Hy@param\Hy@chapapp
    \fi
}{%
    \iftoggle{inappendix}{%true-branch
        \@checkappendixparam{section}%
        \@checkappendixparam{subsection}%
        \@checkappendixparam{subsubsection}%
    }{}%
}{}{\errmessage{failed to patch}}

\newcommand*{\@checkappendixparam}[1]{%
    \def\@checkappendixparamtmp{#1}%
    \ifx\Hy@param\@checkappendixparamtmp
        \let\Hy@param\Hy@appendixstring
    \fi
}
\makeatletter

\newtoggle{inappendix}
\togglefalse{inappendix}

\newif\ifsketchnear
\sketchneartrue
\apptocmd{\appendices}{\toggletrue{inappendix}}{}{\errmessage{failed to patch}}

\newcommand{\refappendix}[1]{\hyperref[#1]{Appendix~\ref*{#1}}}

\usepackage{xspace}
\newcommand{\PIplus}{PI$^{+}$\xspace}

\makeatletter
\patchcmd{\ALG@step}{\addtocounter{ALG@line}{1}}{\refstepcounter{ALG@line}}{}{}
\newcommand{\ALG@lineautorefname}{Line}
\makeatother

\iftoggle{onecolumnversion}{}{
\setlength{\abovedisplayskip}{1ex plus4pt minus2pt}
\setlength{\belowdisplayskip}{\abovedisplayskip}
\setlength{\abovedisplayshortskip}{0pt plus4pt}
\setlength{\belowdisplayshortskip}{1ex plus4pt minus2pt}
}

\usepackage{tikz}
\usepackage{lipsum}

\begin{document}
 
\maketitle
%\copyrightnotice %For copyright compliance

%%%%%%%%%%%%%%%%%%%%%%%%%%%%%%%%%%%%%%%%%%%%%%%%%%%%%%%%%%%%%%%%%%%%%%%%%%%%%%%%
\begin{abstract}
This paper investigates recursive feasibility, recursive robust stability and  near-optimality properties of policy iteration (PI).
For this purpose, we consider deterministic nonlinear discrete-time systems whose inputs are generated by PI for  undiscounted cost functions. 
We first assume that PI is recursively feasible, in the sense that the optimization problems solved at each iteration admit a solution. In this case, we provide novel conditions to establish recursive robust stability properties for a general attractor, meaning that the policies generated at each iteration ensure a robust $\KL$-stability property with respect to a general state measure. We then derive novel explicit bounds on the mismatch between the (suboptimal) value function returned by PI at each iteration and the optimal one. Afterwards, motivated by a counter-example that shows that PI may fail to be recursively feasible, 
we modify PI so that recursive feasibility is guaranteed a priori under mild conditions. This modified algorithm, called \PIplus, is shown to preserve the recursive robust stability when the attractor is compact. Additionally, \PIplus enjoys the same near-optimality properties as its PI counterpart under the same assumptions. Therefore, \PIplus is an attractive tool for generating near-optimal stabilizing control of deterministic discrete-time nonlinear systems.
\end{abstract}

%%%%%%%%%%%%%%%%%%%%%%

\section{Introduction}\label{sec:Intro}

\newcommand{\assname}{Assumption}
\newcommand{\corname}{Corollary}
\newcommand{\thmname}{Theorem}
\newcommand{\propname}{Proposition}
\newcommand{\lemname}{Lemma}
\newcommand{\algorithmname}{Algorithm}
\newcommand{\remname}{Remark}
\newcommand{\defnname}{Definition}
\renewcommand{\sectionautorefname}{Section}
\renewcommand{\subsectionautorefname}{Section}
\renewcommand{\subsubsectionautorefname}{Section}
\makeatletter
\newcommand{\sassname}{SA\@gobble}
\makeatother

Policy iteration (PI) is an optimization algorithm that forms one of the pillars   of   dynamic programming \cite{Bertsekas-book12(adp)}. PI iteratively generates  control laws, also called policies, that converge  to an optimal control law  for general dynamical systems and cost functions under mild conditions, see, e.g.,  \cite{Bertsekas-book12(adp),heydari-acc2016,bian-et-al-aut14,jiang2020learning}. Also, PI may exhibit the attractive feature of converging faster to the optimal value function than its counterpart value iteration (VI) \cite{heydari-acc2016} at the price of more computations. For these reasons, PI attracts a lot of attention both in terms of theoretical investigations see, e.g.,  \cite{heydari-acc2016,BertsekasVIPI,liu-wei-tnnls13,bian-et-al-aut14,modares2013policy,chun2016stability}, and practical applications e.g., \cite{Wang2014cruisePI,Wu2017GasPI,Li2016robotGameTheoryPI,Guo2018PI}.  Nevertheless, several fundamental questions remain largely open  regarding the properties of PI in a control context: (i) its recursive feasibility; (ii) general conditions for recursive robust stability when the attractor is not necessarily a single point but a more general set; (iii) near-optimality guarantees, in particular when the cost function is not discounted. We explain each of these challenges next.  

It is essential that PI is \emph{recursively feasible} in the sense that the optimization problem admits a solution at each iteration. Surprisingly, we have not been able to find general conditions for the recursive feasibility of PI in the literature  when dealing with deterministic nonlinear discrete-time systems with general cost functions, whose state and inputs evolves on a Euclidean space. The only results we came across concentrate on special cases like when the input set is finite \cite{BertsekasVIPI} or the system is linear and the cost is quadratic \cite{Bertsekas-book12(adp)}. The dominant approach in the literature for nonlinear discrete-time systems on Euclidean spaces is thus to assume that the algorithm is recursively feasible, see, e.g., \cite{heydari-acc2016,liu-wei-tnnls13}, or to rely on conditions that are hard to verify a priori in general  as they employ feasibility tests at each iteration  \cite{BertsekasVIPI}. Model predictive control literature recognised a long time ago the importance of recursive feasibility. Hence  we believe that   the recursive feasibility of PI is a property of major importance in view of the burgeoning literature on dynamic programming and reinforcement learning where PI plays a major role \cite{bucsoniu2011approximate,bertsekasbook2019rl}.

A second challenge for PI is related to its application in a control context. In many applications, the closed-loop system must exhibit stability guarantees as: (i) it provides analytical guarantees
on the behavior of the controlled system solutions as time evolves; (ii) it endows the system with  robustness properties
and is thus associated to safety considerations, see, e.g., \cite{berkenkamp2017safe}.   Available results on the stability of systems controlled by PI concentrate on the case where the attractor is a single point, as in, e.g., \cite{bian-et-al-aut14,modares2013policy,chun2016stability,liu-wei-tnnls13}. They exclude set stability, which is inevitable for instance in presence of clock or toggle variables \cite[Examples 3.1-3.2]{Goebel-Sanfelice-Teel-book}, and more generally when the desired operating behaviour of the closed-loop system is given by a set and not a point. Moreover, the commonly used assumptions imposed on the plant model and the stage cost are also subject to some conservatism,  like requiring the stage cost to satisfy positive definiteness properties. In addition, it is essential to ensure that these stability properties are robust, which is not automatically guaranteed, as pointed out in \cite{Grimm-et-al-aut04(examples),Kellett-Teel-siam-jco-05}, and this matter is often eluded in the literature at the exception of the recent work in \cite{Pang2022} in the linear quadratic case. There is therefore a need for  general conditions allowing to conclude robust set stability properties for systems controlled by PI. We further would like these stability properties to be preserved at each iteration; that is, we want to ensure \emph{recursive robust stability}.

Finally, it is important to understand when and how the sequence of value functions  generated by PI at each iteration converges to the optimal value function;  we talk of \emph{near-optimality guarantees}. The literature stands in two ways on this issue. On the one hand, the value functions  are  known to converge monotonically and \emph{point-wisely} to the optimal value function under mild conditions on the model and the cost \cite{bian-et-al-aut14,heydari-acc2016}; uniform convergence properties are only ensured, as far as we know, for discounted costs functions \cite{Bertsekas-book12(adp)}.   On the other hand, it is important to be able to evaluate the mismatch between the returned value function at each iteration and the optimal one. These computable or explicit bounds on the mismatch of generated value functions are vital to  decide when to stop iterating the algorithm. Existing results  concentrate on discounted costs \cite{Bertsekas-book12(adp)}, which are not always natural in control applications. In the discounted setting, the provided near-optimality bounds explode when the discount factor $\gamma\in(0,1)$ converges to one \cite{Bertsekas-book12(adp)}, while to ensure stability, in opposition, $\gamma$ should be close enough to one in view of \cite{postoyan_stability_discounted,granzotto2020finite}.  Hence, there is a need for stronger near-optimality guarantees for PI. In particular, we seek to develop results  that provide computable near-optimality bounds, without relying on a discount factor, and to provide conditions under which the sequence of constructed value functions satisfies a \emph{uniform} monotonic convergence property towards the optimal one.

In this context, we consider deterministic nonlinear discrete-time systems, whose inputs are generated by PI for an undiscounted infinite-horizon cost function. We first assume that PI is recursively feasible and we provide general conditions inspired from the model predictive control literature \cite{Grimm-et-al-tac2005} to ensure the recursive robust stability of the closed-loop system at each iteration, where the attractor is a set. These conditions relate to the detectability of the system with respect to the stage cost and the stabilizing property of the initial policy. We then exploit these stability properties to derive explicit and computable bounds on the mismatch between the optimal value function and the value function obtained by PI at each iteration. We also show that the sequence of value functions satisfies a uniform convergence property towards the optimal value function by exploiting stability. 

Afterwards, we show via a counter-example that PI may actually fail to be recursively feasible under commonly used assumptions. We thus propose to modify the original formulation of PI so that we can guarantee the recursive feasibility of the algorithm under mild conditions on the  model, the stage cost and the input set. We call this new algorithm PI plus (\PIplus). \PIplus differs from PI in two aspects.  First, an (outer semicontinuous) regularization is performed at the so-called improvement step, which is a common technique in the discontinuous/hybrid systems literature \cite{Goebel-Sanfelice-Teel-book}. Second, instead of letting  the algorithm select any policy that minimizes the (regularized) improvement step, we select any of those generating the smallest cost that we aim to minimize thereby requiring an extra layer of computation. In this paper, we do not address the question of the practical implementation of \PIplus, which is left for future work. Instead, we concentrate on the methodological challenges raised by the algorithm.  We then prove that \PIplus is indeed recursively feasible, and that it  preserves the recursive robust stability when the attractor is compact as well as the  near-optimality properties established for PI. Compared to our preliminary work in  \cite{granzotto-et-all-PI-CDC}, novel elements include the results for PI, the robust stability analysis, the fact that the admissible input set can be state-dependent, and new technical developments to derive less conservative near-optimality bounds. Moreover, the full proofs are provided.

The rest of the paper is organized as follows.  Preliminaries  are given in \autoref{sec:notation}. The analysis of PI is carried out in \autoref{sec:sys-cost-problem}. The new algorithm \PIplus and its properties are  presented in \autoref{sec:mainresults}. We defer the robustness analysis of the stability properties ensured by PI and \PIplus in \autoref{sec:robust}.  Concluding remarks are provided in \autoref{sec:conclusion}. In order to streamline the presentation, the proofs are given in the appendices. 

%%%%%%%%%%%%%%%%%%%%%%

\section{Preliminaries} \label{sec:notation}

In this section, we define the  notation,  provide important definitions and  formalize the problem.

\subsection{Notation} \label{subsec:notation}
Let $\R:=(-\infty,\infty)$,  
$\Rlo:=[0,\infty)$,   $\overline{\R}:=[-\infty,\infty]$,  $\Zo:=\{0,1,2,\ldots\}$ and  $\Zp:=\{1,2,\ldots\}$. The notation $(x,y)$ stands for $[x^{\top},\,y^{\top}]^{\top}$, where $x\in\Rl{n}{}$,  $y\in\Rl{m}{}$ and $n,m\in\Zp$. The Euclidean norm of a vector $x\in\Rl{n}{}$ with $n\in\Zp$ is denoted by $|x|$ and the distance of $x\in\Rl{n}{}$ to a non-empty set $\mathcal{A}\subseteq\Rl{n}{}$ is denoted by  $|x|_{\mathcal{A}}:=\inf\{|x-y|\,:\,y\in\mathcal{A}\}$. 
The unit closed ball of $\R^{n}$ for $n\in\Zp$ centered at the origin is denoted by $\mathbb{B}$.
We consider $\K$, $\K_\infty$ and $\KL$ functions as defined in \cite[Section 3.5]{Goebel-Sanfelice-Teel-book}. We write $\beta\in{\exp}{-}\KL$ when $\beta(s_1,s_2)=\lambda_1 s_1 e^{-\lambda_2 s_2}$ for some $\lambda_1\in[1,\infty)$ and $\lambda_2>0$ for any $(s_1,s_2)\in\Rlo^2$. For any set $\mathcal{A}\subseteq\R^{n}$, $x\in\R^{n}$, the \emph{indicator function} $\delta_{\mathcal{A}}:\R^{n}\to\overline{\R}$ is defined as $\delta_{\mathcal{A}}(x)=0$ when $x\in\mathcal{A}$ and $\delta_{\mathcal{A}}(x)=\infty$ when $x\notin\mathcal{A}$ as in \cite{Rockafellar-Wets-book}. Moreover, when $\mathcal{A}$ is closed, we say $\sigma:\R^{n}\to\Rlo$ is a \emph{proper indicator of set $\mathcal{A}$} whenever  $\sigma$ is continuous and there exist $\underline{\alpha}_\sigma,\overline{\alpha}_\sigma\in\K_\infty$ such that $\underline{\alpha}_\sigma(|x|_{\mathcal{A}})\leq\sigma(x)\leq\overline{\alpha}_\sigma(|x|_{\mathcal{A}})$ for any $x\in\R^{n}$.
The identity map  from $\Rlo$ to $\Rlo$ is denoted by $\mathbb{I}$, and the zero map from $\Rlo$ to $0$ by $\bm{0}$. Let $f: \Rlo \to\Rlo$. We use $f^{(k)}$ for the composition of function $f$
to itself $k$ times, where $k\in\Zo$  and $f^{(0)}:=\mathbb{I}$. Given a set-valued map $S:\R^n \rightrightarrows\R^m$, a selection of $S$ is a single-valued mapping $s:\dom S \to \R^m$ such that $s(x)\in S(x)$ for any $x\in \dom S$. For the sake of convenience, we write $s\in S$ to denote a selection $s$ of $S$. We also employ the following  definition from \cite[Def. 1.16]{Rockafellar-Wets-book}.

\begin{defn}[uniform level boundedness] \label{def:levelboundeduniform}  A function $f:\R^{n}\times\R^{m}\to \overline{\R}$ where $n,m\in\Zp$ with values $f(x,u)$ is \emph{level-bounded in $u$, locally uniform in $x$} if for each $x\in\R^{n}$ and $\alpha\in\R$  there is a neighborhood   $\mathcal{S}$ of $x$ along a bounded set $B\subset\R^{n}$ such that $\{u\in\R^{m}\,:\,f(z,u)\leq\alpha\}\subset B$ for any $z\in \mathcal{S}$. 
$\hfill\Box$
\end{defn}
\subsection{Plant model and cost function} \label{subsec:optimalcontrolproblem}
Consider the plant model
\begin{equation}
x(k+1) =  f(x(k),u(k)),
\label{eq:sys-plant}
\end{equation}
where $x\in\R^{n_x}$ is the state, $u\in\mathcal{U}(x)\subseteq\R^{n_u}$ is the control input, the time-step is $k\in\Zo$, $\mathcal{U}(x)$ is a non-empty set of admissible inputs for state $x\in\R^{n_x}$, and $n_x,n_u\in\Zp$. We wish to find, for any given $x\in\R^{n_x}$, an   infinite-length sequence of admissible inputs $\bm{u}=(u(0),u(1),\dots)$, that minimizes the   infinite-horizon   cost
\begin{equation}
J(x,\bm{u})  :=   \sum_{k=0}^{\infty}\ell(\phi(k,x,\bm{u}|_{k}),u(k)),
\label{eq:J}
\end{equation} 
where {$\ell:\R^{n_x}\times\R^{n_u}\to\Rlo$} is a non-negative stage cost and  $\phi(k,x,\bm{u}|_{k})$ is the solution to \eqref{eq:sys-plant} at the {$k^{\text{th}}$-step}, initialized at $x$ at time 0  with inputs $\bm{u}|_k:=(u(0),\ldots,u({k-1}))$.   The   minimum of $J(x,\cdot)$ is denoted as  
\begin{equation} \label{eq:Vstar}
    V^\star(x):= \min_{\bm{u}} J(x,\bm{u})
\end{equation}
for any $x\in\R^{n_x}$,  where $V^\star$ is the    \emph{optimal value function} associated to the minimization of \eqref{eq:J}.
\begin{sass}[\autoref{SA:well-posed}]\label{SA:well-posed} 
For any $x\in\R^{n_x}$, there exists an optimal sequence of admissible inputs $\bm{u}^\star(x)$ such that $V^{\star}(x)=J(x,\bm{u}^\star(x))< {+\infty}$ and for any infinite-length sequence of admissible inputs $\bm{u}$, $V^{\star}(x)\leq J(x,\bm{u})$.\mbox{}\hfill $\Box$
\end{sass}

Conditions to ensure \autoref{SA:well-posed} can be found in, e.g., \cite{Keerthi-Gilbert-tac85}. Given \eqref{eq:Vstar}, we  define the set of optimal inputs as
\begin{equation} \label{eq:Hstar}
    H^\star(x):= \argmin_{u\in\mathcal{U}(x)}\left\{\ell(x,u)+ V^\star(f(x,u))\right\}.
\end{equation}

To compute $H^\star$ in \eqref{eq:Hstar} for the general dynamics in \eqref{eq:sys-plant} is notoriously hard. Dynamic programming  provides algorithms to iteratively  obtain feedback laws, which instead converge to $H^{\star}$ \cite{bertsekasbook2019rl}. A fundamental algorithm of dynamic programming is PI, which is presented and analysed in the next section. Before that, we introduce some notation, which will be convenient in the sequel. Given a feedback law $h:  \R^{n_x}\to\R^{n_u}$ that is admissible, i.e., $h\in\mathcal{U}$,  we denote the solution to system \eqref{eq:sys-plant} in closed-loop with feedback law $h$ at time $k\in\Zo$ with initial condition $x$ at time 0 as $\phi(k,x,h)$. Likewise, $J(x,h)$ is the cost induced by $h$ at initial state $x$, i.e., $J(x,h)=\sum_{k=0}^\infty \ell(\phi(k,x,h),h(\phi(k,x,h)))$. In this way, we have that for any selection $h^\star\in H^\star$, $V^\star(\cdot)=J(\cdot,h^\star)$, as $H^\star(x)$ is non-empty for any $x\in\R^{n_x}$ by \autoref{SA:well-posed}.

%%%%%%%%%%%%%%%%%%%%%%
\section{Policy Iteration}
\label{sec:sys-cost-problem}

We recall in this section  the original formulation of PI and we assume the algorithm is feasible at any iteration. We then establish  novel recursive stability and near-optimality guarantees, assuming a detectability property is satisfied and the initial policy is stabilizing. Finally, we present an example where PI is not recursively feasible despite supposedly favorable properties, thereby motivating the need to modify PI to overcome this issue, which will be the topic of \autoref{sec:mainresults}. As mentioned in the introduction, the robustness of the stability properties is analyzed in \autoref{sec:robust}.

\subsection{The algorithm}
\label{subsec:piClassical}

PI is presented  in \autoref{algo:PI}. Given an initial admissible policy $h^0$, 
PI generates at each iteration $i\in\Zo$  a policy $h^{i+1}$ with cost $V^{i+1}(x):=J(x,h^{i+1})$ and it can be proved that $V^{i+1}(x) \leq V^{i}(x)$ for all $x\in\R^{n_x}$ \cite[Section 4.2]{Sutton-Barto-book2018}. 
This is done via the  improvement step in \eqref{eq:PI-ith-improvement-classic}.
The policy $h^{i+1}$ obtained at  iteration $i+1$ is an arbitrary selection of $H^{i+1}$ in  \eqref{eq:PI-ith-improvement-classic} where $H^{i+1}$  may be set-valued.  
We then evaluate the cost induced by $h^{i+1}$, namely $V^{i+1}=J(\cdot,h^{i+1})$, this is the evaluation step in \eqref{eq:PI-ith-evaluation-classic}.
By doing so repeatedly, $V^i$ converges to the optimal value function $V^\infty=V^\star$ under mild conditions,  see \cite{Bertsekas-book12(adp)}. 

\begin{algorithm}[t]
\renewcommand{\hypcapspace}{2\baselineskip}
\capstart
\renewcommand{\hypcapspace}{0.5\baselineskip}
 \caption{\label{algo:PI} Policy Iteration (PI)}
 \begin{algorithmic}[1]
 \renewcommand{\algorithmicrequire}{\textbf{Input:}}
 \renewcommand{\algorithmicensure}{\textbf{Output:}}
 \REQUIRE $f$ in \eqref{eq:sys-plant}, $\ell$ in \eqref{eq:J}, initial policy $h^0\in\mathcal{U}$
 \ENSURE  Policy $h^\infty$, cost $V^\infty$
 \STATE \textbf{Initial evaluation step:} for all $x\in\R^{n_x}$, 
  \begin{equation}V^0(x) := J(x,h^0).\tag{PI.1} \label{eq:PI-0th-evaluation}
 \end{equation}
  \FOR{$i\in\Zo$}
  \STATE \textbf{Policy improvement step:}  for all $x\in\R^{n_x}$,
  \begin{equation}
     H^{i+1}(x) := \argmin\limits_{u\in\mathcal{U}(x)} \{\ell(x,u)+V^{i}(f(x,u))\}.
 \tag{PI.2} \label{eq:PI-ith-improvement-classic}
\end{equation}
\STATE  \textbf{Select} $h^{i+1}\in H^{i+1}$. \label{algo-line:PI-policy-choice}
\STATE \textbf{Policy evaluation step:} for all $x\in\R^{n_x}$,
 \begin{equation}
V^{i+1}(x):= J(x,h^{i+1}).  
\tag{PI.3}
\label{eq:PI-ith-evaluation-classic}
\end{equation}
\ENDFOR
 \STATE \textbf{return} $h^\infty \in H^\infty$ and $V^\infty$.
 \end{algorithmic}
\end{algorithm}

\begin{rem}\label{rem:cstop}
In practice, PI is often stopped at some iteration, typically by  looking at the difference between $V^i$ and $V^{i-1}$ for some $i\in\Zp$. We will return to this point in \autoref{rem:cstop-PI} in \autoref{ssubsec:near-opti-PI}. \mbox{}\hfill$\Box$
\end{rem}

\subsection{Desired properties} \label{ss:desiredproperties}
For the remainder of this section, we proceed as is often done in the literature and assume \autoref{algo:PI} is recursively feasible, see, e.g., \cite{heydari-acc2016,liu-wei-tnnls13}, in the sense that the optimization  problem in \eqref{eq:PI-ith-improvement-classic} admits a solution for any $x\in\R^{n_x}$ at any iteration $i\in\Zp$, which is formalized below.

\begin{ass}
Set-valued map $H^i(x)$ is non-empty for any $i\in\Zp$ and  $x\in\R^{n_x}$. \mbox{}\hfill$\Box$\label{ass:rfeasibility}
\end{ass}

We note that verifying \autoref{ass:rfeasibility} is hard in general. At a given iteration $i+1$ with $i\in\Zo$,  a sufficient condition for $H^{i+1}(x)$  to be non-empty for any $x\in\R^{n_x}$ is establishing the lower semicontinuity of map $(x,u)\mapsto \ell(x,u)+V^i(f(x,u))+\delta_{\mathcal{U}(x)}(u)$ on $\R^{n_x}\times\R^{n_u}$, see, e.g., \cite{BertsekasVIPI}. However, this  is hard to check in advance as $V^i$ is not known a priori. We will return to the question of the recursive feasibility of PI in \autoref{ssec:counterexample}.

Under \autoref{ass:rfeasibility}, the goal of this section is to establish the recursive stability and near-optimality bounds for PI, in particular we aim at showing that PI is 
\begin{itemize}
    \item \emph{recursively stabilizing}, i.e., if $h^{0}$ stabilizes system~\eqref{eq:sys-plant}, then this property is also ensured by any $h^{i}\in H^{i}$ for any $i\in\Zp$, in the sense that  the difference inclusion%
\begin{equation}
    x(k+1)\in  f(x(k),H^{i}(x(k))) =: F^{i}(x(k)) \label{eq:auto-sys-PI}
\end{equation}
exhibits desirable set stability properties;
    \item \emph{near-optimal} in the sense that we have guaranteed bounds  on $V^{i}-V^{\star}$ for any $i\in\Zo$, despite the fact that $V^{\star}$ in \eqref{eq:Vstar} is typically unknown.
\end{itemize}

These results rely on assumptions given in the next section. For convenience,  solutions to system \eqref{eq:auto-sys-PI} are denoted in the sequel as   $\phi^{i}(\cdot,x)$ when initialized at some $x\in\R^{n_x}$ for any $i\in\Zo$.

\subsection{Standing assumptions} \label{subsec:assumptions}

To define stability, we use a continuous function $\sigma:\R^{n_x}\to\Rlo$ that  serves as a state ``measure''  relating the distance of the state to a given attractor where $\sigma$ vanishes.  By stability, we  mean that there exists $\beta\in\KL$ (independent of $i$) such that, for any $x\in\R^{n_x}$, $i\in\Zo$, any solution $\phi^i$ to \eqref{eq:auto-sys-PI} verifies, for any $k\in\Zo$,
\begin{equation}\label{eq:PIplusstab-0th}
        \sigma(\phi^{i}(k,x))\leq\beta(\sigma(x),k).
    \end{equation}
Property \eqref{eq:PIplusstab-0th} is a $\KL$-stability property of system \eqref{eq:auto-sys-PI} with respect to $\sigma$. When $\sigma$ is a proper indicator function of a closed set $\mathcal{A}\subseteq\R^{n_x}$,  the uniform global asymptotic stability of  set  $\mathcal{A}=\{x\in\R^{n_x}\,:\,\sigma(x)=0\}$  is guaranteed by \eqref{eq:PIplusstab-0th}. When $\sigma(x)=|x|^{a}$ for any $x\in\R^{n_x}$ with $a\geq 1$, the uniform global asymptotic stability of the origin $x=0$ is ensured by \eqref{eq:PIplusstab-0th}, for instance. Function $\sigma$ is thus convenient to address stability properties for general attractors.

We make the next detectability assumption on system~\eqref{eq:sys-plant} and stage cost $\ell$ consistently with e.g., \cite{Grimm-et-al-tac2005,postoyan_stability_discounted,granzotto2020finite,kalman_contributions_1960}.
\begin{sass}[\autoref{SA:detect-control}] \label{SA:detect-control} There exist  a  continuous function $W:\R^{n_x}\to\Rlo$, $\alpha_W, \chi_W \in \Kinf$ and $\overline\alpha_W:\Rlo\to\Rlo$ continuous, nondecreasing and zero at zero, such that, for~any $(x,u)\in \mathcal{W}:=\{(x,u)\in\R^{n}\times\R^{m}\,:\,u\in\mathcal{U}(x)\}$,\\
\begin{minipage}{\columnwidth}
\begin{equation}
\begin{split}
W(x) &  {} \leq  \overline{\alpha}_W(\sigma(x))\\
W(f(x,u)) - W(x) & {} \leq  -\alpha_W(\sigma(x)) +  \chi_W(\ell(x,u)).
\end{split}\label{eq:a-detectability}
\end{equation}
\mbox{} \hfill$\Box$
\end{minipage}
\end{sass}

\autoref{SA:detect-control} is a   detectability property of system \eqref{eq:sys-plant} and stage cost $\ell$, see \cite{Grimm-et-al-tac2005,Grune_diss_detec_2019} for more details. \autoref{SA:detect-control} holds for instance with $W=0$ when  $\{x\in\R^{n_x}\,:\,\sigma(x)=0\}$ is compact, $\ell(x,u) =  \ell_1(x)+\ell_2(x,u)$ with $\ell_1$ continuous, positive definite with respect to the set $\{x\in\R^{n_x} : \sigma(x)=0\}$, $\ell_1(x)\to \infty$ as $|x|\to\infty$, and $\ell_2(x,u)\geq 0$ for any $(x,u)\in\mathcal{W}$. Note that \autoref{SA:detect-control} relaxes the requirement that $\ell(\cdot,u)$ be positive definite for any $u\neq 0$ as found in, e.g., \cite{bian-et-al-aut14,modares2013policy,chun2016stability,liu-wei-tnnls13}, and $\ell$ is not required to satisfy convexity properties.

Finally, like in e.g., \cite{heydari-acc2016,BertsekasVIPI,liu-wei-tnnls13}, we assume  that we  initialize the algorithm  with a stabilizing feedback law $h^0$. In particular, we make the next assumption. 

\begin{sass}[\autoref{SA:stabilizing-initial-policy}]\label{SA:stabilizing-initial-policy} There exists $\overline{\alpha}_V\in\Kinf$ such that, for any $x\in\R^{n_x}$, $V^{0}(x)=J(x,h^0)\leq\overline{\alpha}_V(\sigma(x))$.
\mbox{} \hfill $\Box$

\end{sass}

\autoref{SA:detect-control} and \autoref{SA:stabilizing-initial-policy} are related to the stability property of system~\eqref{eq:sys-plant} in closed-loop with $h^0$ as in \eqref{eq:PIplusstab-0th}. This is established below. 

\subsection{Results}\label{ssec:piprop}

We now establish the desired properties listed in \autoref{ssec:piprop}. The proofs of the forthcoming results follow very similar lines as the equivalent ones stated in \autoref{sec:mainresults} for \PIplus. For this reason, the proofs are carried out in details for the results of \autoref{sec:mainresults} in \hyperref[appendix:recursive]{Appendices \ref*{appendix:recursive}}  \hyperref[appendix:secnear-opti]{and \ref*{appendix:secnear-opti}}, and their application to the results of the present section is discussed in \refappendix{appendix:PI}.

\subsubsection{Recursive stability}
The next theorem establishes recursive stability.

\begin{thm}\label{thm:rstability-PI}
Suppose \autoref{ass:rfeasibility} holds. For any $i\in\Zo$, for any $x\in\R^{n_x}$, solution $\phi^i$ to \eqref{eq:auto-sys-PI} and $k\in\Zo$,
\begin{samepage}
\begin{equation}\label{eq:PIplusstab-PI}
        \sigma(\phi^{i}(k,x))\leq\beta(\sigma(x),k),
    \end{equation}
    where\footnote{\label{foot:non-decreasing}To simplify presentation, we assume $\mathbb{I}-\widetilde{\alpha}_Y$ is non-decreasing. If this is not the case, we can always replace $(\mathbb{I}-\widetilde{\alpha}_Y)^{(k)}$ by $s\mapsto\max_{\hat{s}\in[0,s]}(\mathbb{I}-\widetilde{\alpha}_Y)^{(k)}(\hat{s})$ in the expression of $\beta$.} $\beta:(k,s)\mapsto\underline{\alpha}_Y^{-1} \circ (\mathbb{I}-\widetilde{\alpha}_Y)^{(k)}\circ \overline{\alpha}_Y(s)\in\KL$ with $\underline{\alpha}_Y,\overline{\alpha}_Y,\alpha_Y$ in Table \ref{table:functions-theorems}.  \mbox{}\hfill $\Box$
\end{samepage} 
\end{thm}
\autoref{thm:rstability-PI} ensures the desired $\KL$-stability property of system~\eqref{eq:auto-sys-PI} with respect to $\sigma$ at any iteration $i\in\Zo$, and thus at $i=0$, which confirms that $h^0$ is stabilizing as mentioned above. It is important to note that $\beta$ in \eqref{eq:PIplusstab-PI}  is independent of the number of iterations $i$, which makes the stability property uniform with respect to $i$.

Under extra conditions, we can derive an exponential stability result.

\begin{cor} \label{cor:expresult-PI}
Suppose \autoref{ass:rfeasibility} holds and that there exist  $c_W,a_W,\overline a_V>0$ and $\overline a_W\geq0$ such that $\chi_W(s)\leq c_W s$, $\alpha_W(s)\geq a_W s$, $\overline\alpha_V(s)\leq \overline a_V s$, $\overline\alpha_W(s)\leq \overline a_W s$ for any $s\geq 0$, where $\chi_W,\alpha_W,\overline\alpha_W$ come from \autoref{SA:detect-control} and $\overline\alpha_V$ comes from \autoref{SA:stabilizing-initial-policy}. Then,
\autoref{thm:rstability-PI} holds with $\beta:(s,k)\mapsto\tfrac{\overline a_Y}{\underline a_Y}(1-\widetilde{a}_Y)^k s\in{\exp}{-}\KL$, and $\widetilde{a}_Y,\overline a_Y,\underline a_Y$ in Table \ref{table:functions-theorems}.
\hfill $\Box$
\end{cor}

\subsubsection{Near-optimality properties}\label{ssubsec:near-opti-PI}

Next, we establish near-optimality properties of PI.

\begin{thm} \label{thm:near-opti-PI} Suppose \autoref{ass:rfeasibility} holds. For any $i\in\Zo$ and $x\in\R^{n_x}$, 
\begin{equation} \label{eq:abstract-near-pi}
\begin{array}{rlll}
(V^i-V^\star)(x) & \leq & (V^0-V^\star)(\phi(i,x,h^\star)),
\end{array}
\end{equation}
     where  $h^\star \in H^\star$ from \eqref{eq:Hstar}. Moreover, for any $i\in\Zo$ and $x\in\R^{n_x}$,
     \begin{equation}
     \begin{array}{rlll}
         V^i(x)-V^\star(x) & \leq & \widetilde{\alpha}\left(\beta(\sigma(x),i)\right), \label{eq:explicit-near-PI}
     \end{array}
     \end{equation}
     where $\beta\in\KL$ comes from \autoref{thm:rstability-PI} and $\widetilde{\alpha}=\max_{\hat{s}\in[0,s]} \widehat{\alpha}(\hat{s})$ is non-decreasing and positive definite function with  $\widehat{\alpha}$ in Table \ref{table:functions-theorems}. \hfill $\Box$
\end{thm}

\autoref{thm:near-opti-PI} provides novel characterisations of the near-optimality properties of PI. In  \eqref{eq:abstract-near-pi}, $(V^i-V^\star)(x)$ is the near-optimality error term at iteration $i\in\Zo$ and state $x\in\R^{n_x}$. This error is upper-bounded in \eqref{eq:abstract-near-pi} by $(V^0-V^\star)(\phi(i,x,h^\star))$, which is the  ``initial'' near-optimality error term for $i=0$, but evaluated at state $\phi(i,x,h^\star)$ instead of $x$. In its turn, state $\phi(i,x,h^\star)$ corresponds  to the $i^\text{th}$ time-step of the solution of \eqref{eq:sys-plant} initialized at $x$ in closed-loop with an \emph{optimal} (and typically unknown) policy $h^\star\in H^\star$. This bound decreases to zero point-wisely as $i\to\infty$ thanks for the stability property of system \eqref{eq:sys-plant} in closed-loop with \eqref{eq:Hstar}, which follows similarly as \autoref{thm:rstability-PI}.  However, the upper-bound in \eqref{eq:abstract-near-pi} is typically unknown. To overcome
this possible issue, a  conservative upper-bound of $(V^0-V^\star)(\phi(i,x,h^\star))$ in the form of $\widetilde{\alpha}\left(\beta(\sigma(x),i)\right)$ is given in \eqref{eq:explicit-near-PI}, where $\beta$ and $\widetilde{\alpha}$  come from \hyperref[thm:rstability-PI]{Theorems \ref*{thm:rstability-PI}} \hyperref[thm:near-opti-PI]{and \ref*{thm:near-opti-PI}} and can be both computed. The fact that function $\beta$ in \eqref{eq:PIplusstab-PI}  is independent of the number of iterations $i$ is vital for \eqref{eq:explicit-near-PI}.  Indeed, as a result, the upper-bound is ensured to converge to zero as $i$ increases to infinity. We emphasize that the above near-optimality properties exploit stability to provide explicit bounds as in \eqref{eq:explicit-near-PI}. This is in contrast to the literature, which relies on discount factors to derive contractive properties of the sequence of value functions generated by PI \cite{Bertsekas-book12(adp)}.

\begin{rem}\label{rem:cstop-PI}
The bound in \eqref{eq:explicit-near-PI} can be used to design stopping criteria for PI for a given  near-optimality guarantee. To see this, consider a desired near-optimality target $\varepsilon_\text{target}:\Rlo\to\Rlo$,  which may only vanish at the origin. It suffices to iterate PI until $i\geq i^\star$, where $i^\star$ is such that
\begin{equation}
    \widetilde{\alpha}\left(\beta(\sigma(x),i^\star)\right)\leq  \varepsilon_\text{target}(\sigma(x))
\end{equation} 
for every $x\in\R^{n_x}$. As a result, for any $i\geq i^\star$ and $x\in\R^{n_x}$, it holds that $V^i(x)-V^\star(x)\leq \varepsilon_\text{target}(\sigma(x))$ as desired. Moreover, $\varepsilon_\text{target}(\sigma(x))$ may be as small as desired for $x$ in a given compact set, provided $i^\star$ is sufficiently large. \mbox{}\hfill$\Box$ 
\end{rem}

A direct consequence of \autoref{thm:near-opti-PI} is that the sequence of cost functions $V^i$  satisfies a uniform convergence property towards $V^\star$, as formalized next.

\begin{prop} \label{prop:uniform-PI}
Suppose \autoref{ass:rfeasibility} holds. The sequence of functions $V^i$ monotonically uniformly converges to $V^\star$ on level-sets  of $\sigma$, i.e.,:
\begin{itemize}
    \item[(i)] for all $x\in\R^{n_x}$, $V^{i+1}(x)\leq V^{i}(x)$;
    \item[(ii)] for any $\Delta,\delta>0$, there exists $i^\star\in\Zp$ such that, for all $x\in\{z\in\R^{n_x}\,:\, \sigma(z)\leq\Delta\}$ and $i\geq i^\star$, $V^i(x)-V^\star(x) \leq  \delta$. 
$\mbox{}\hfill\Box$
\end{itemize}
\end{prop}

\autoref{prop:uniform-PI} implies a monotonic uniform convergence property of $V^i$ to $V^\star$.  This is an additional benefit of our analysis compared to the existing PI literature, for which only monotonic point-wise convergence of the value functions to the optimal one is guaranteed in general \cite{BertsekasVIPI,heydari-acc2016}  when the cost is not discounted as in \eqref{eq:J}.

The results of this section rely on \autoref{ass:rfeasibility}, namely that PI is recursively feasible. We show next that \autoref{ass:rfeasibility} may fail to hold even when the system and the cost satisfy supposedly favorable properties. 

\begin{table}[ht] \vspace{-2ex}
\begin{tabularx}{\columnwidth}{rllX}
General case \\
\midrule
$\widetilde{\alpha}_Y\,\,\,\,=$ & $\alpha_Y \circ \overline{\alpha}_Y^{-1}$\\
$\alpha_Y\,\,\,\,=$ & $q_{W}(\frac{1}{4}\alpha_W)\frac{1}{4}\alpha_W$\\
$\overline\alpha_Y\,\,\,\,=$ & $\rho_V\circ\overline\alpha_V+\rho_W\circ\overline\alpha_W$\\
$\underline\alpha_Y\,\,\,\,=$ &$\min\left\{\rho_V\circ\chi_W^{-1}(\frac{1}{2}\alpha_W),\rho_W(\frac{1}{2}\alpha_W)\right\}$ \\
$q_V\,\,\,\,=$ & $\dst 2\chi_W(2\mathbb{I})$ \\
$q_W\,\,\,\,=$ & $\frac{1}{2}[\chi_W+(\overline\alpha_W+\mathbb{I})\circ\alpha_W^{-1}(2\chi_W)]^{-1}$ \\
$\rho_V(s)\,\,\,\,=$ & $\int_{0}^{s}q_V(\tau)d\tau$ $\,\,\,\,\forall s\geq 0$ \\
$\rho_W(s)\,\,\,\,=$ & $\int_{0}^{s}q_W(\tau)d\tau$  
 $\,\,\,\forall s\geq 0$\\
$\widehat{\alpha}\,\,\,\,=$ & $\overline{\alpha}_{V}-\chi_W^{-1}\left(\max\{0,\alpha_W-\overline{\alpha}_W\}\right)$\\ 
 \midrule
\end{tabularx}
\begin{tabularx}{\columnwidth}{rlrl}
\multicolumn{2}{l}{When $\chi_W\leq \mathbb{I}$} & \multicolumn{2}{l}{Under the  conditions of \autoref{cor:expresult-PI}}  \\
\midrule
$\widetilde{\alpha}_Y\,\,\,\,=$ & $\alpha_Y \circ \overline{\alpha}_Y^{-1}$ &
$\widetilde{a}_Y\,\,\,\,=$ & $ \nicefrac{a_Y}{\overline{a}_Y}$\\
$\alpha_Y\,\,\,\,=$ & $\alpha_W$ &
$a_Y\,\,\,\,=$ & $a_W$\\ 
$\overline{\alpha}_Y\,\,\,\,=$ & $\overline{\alpha}_V+\overline{\alpha}_W$ &  $\overline{a}_Y\,\,\,\,=$ & $\overline{a}_V+\overline{a}_W$\\ 
$\underline{\alpha}_Y\,\,\,\,=$ & $\alpha_W$ &  $\underline{a}_Y\,\,\,\,=$ & $a_W$ \\
$\widehat{\alpha}\,\,\,\,\,\,=$ & $\min\{\overline{\alpha}_V,\overline{\alpha}_Y-\underline{\alpha}_Y\}$ & $\widehat{\alpha}\,\,\,\,\,\,=$ & $\min\{\overline{a}_V,\overline{a}_Y-\underline{a}_Y\} \cdot \mathbb{I}$ \\
\midrule
\end{tabularx}
\caption{\label{table:functions-theorems}Functions in \hyperref[thm:rstability-PI]{Theorems \ref*{thm:rstability-PI}} \hyperref[thm:near-opti-PI]{and \ref*{thm:near-opti-PI}}.}
\end{table}

\subsection{Recursive feasibility: a counter-example for PI}\label{ssec:counterexample}
Consider the input-affine system
\begin{equation} \label{eq:syscounter}
    x(k+1)= (1-u(k))\max\{0,|x(k)|-1\},
\end{equation}
with $x\in\R$, $u\in\mathcal{U}(x)=[-\delta,1]$ with $\delta=\frac{1}{100}$. Notice that $f:(x,u)\mapsto (1-u(k))\max\{0,|x(k)|-1\}$ is continuous on $\R^{n_x}\times\R^{n_u}$ and that the admissible set of input is compact and convex. Stage cost $\ell$ is  defined as
    $\ell(x,u) = 3|x|g_1(u)+\left(|x|+\frac{7}{4}|x|^2\right)g_2(u),$
where $g_1(u) := \max\{\min\{2(1-u),1\},0\}$ and $g_2(u) := \max\{\min\{2u,1\},0\}$, for any $x\in\R$ and $u\in\mathcal{U}(x)$.  Note that $\ell(x,u)\geq0$ for any $x\in\R$, $u\in\mathcal{U}(x)$ and $\ell(0,u)=0$  for any $u\in\mathcal{U}(0)$. 

Let $h^0(x)=0$ for all $x\in\R$. We obtain from \autoref{algo:PI} that $V^0(x)=3|x|$ for $|x|\in[0,1]$, $V^0(x)=6|x|-3$ for $|x|\in[1,2]$, $V^0(x)=9|x|-9$ for $|x|\in[2,3]$ and so on,  hence $V^0$ is continuous.  As a result, since $\mathcal{U}(x)$ is compact, $H^1(x)$ in \eqref{eq:PI-ith-improvement-classic} is non-empty for any $x\in\R$. Consider $\bar{x}=\frac{18}{7}$, we have that $H^1(\bar{x})=\{0,1\}$:  $H^1$ is set-valued at $\bar{x}$. This implies that we can consider two distinct   policies $h^1,h^{1'}\in H^1$ such that
$h^1(\bar{x})=1$ and $h^{1'} (\bar{x})=0$.     For $x\in(\bar{x},\infty)$  $H^1(x)=0$, for  $x\in(0,\bar{x})\setminus\{2\}$  $H^1(x)=1$, $H(2)=\{0,1\}$, $H^1(0){=}\mathcal{U}(0)$ and $H(x)=H(-x)$ for $x<0$. %V^0(x)= K(x)*x-m(x), where K(n)=3*(n+1), m(n)=3/2  n (n+1) for n non-negative integer, hence V^0(n) = 3/2  n (n+1) which allows for V^0(n)<ell(n,1) for n>2.

\begin{figure}
\centering 
    \includegraphics[width=0.7\columnwidth,trim={0 0.1cm 0cm 2.8cm},clip]{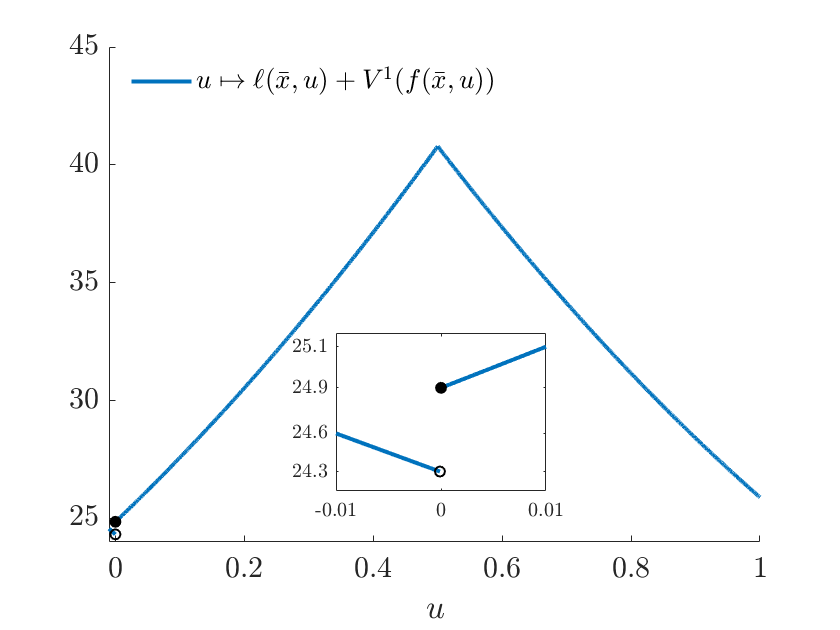}
    
    \caption{$u\mapsto\ell(\bar x+1, u) + V^{1}(f(\bar x+1,u))$  has no minimum.}
    \label{fig:unfeasible}
\end{figure}

Let $V^1(x) := J(x,h^1)$ and $V^{1'}(x) := J(x,h^{1'})$. We note that $V^1(\bar x) = \ell(\bar{x},1)=\frac{396}{28}>\frac{381}{28}=\ell(\bar{x},0)+\ell(\bar{x}-1,1)= V^{1'}(\bar x)$. Hence the policies $h_1$ and $h_1'$ lead to different value functions. Consider $h^1$ and $V^1$.  We see in \autoref{fig:unfeasible} that $u\mapsto\ell(\bar x+1, u) + V^{1}(f(\bar x+1,u))$  has no minimum over $\mathcal{U}(\bar x+1)$, but only an infimum at $u=0$. As a result, the minimization step in \eqref{eq:PI-ith-improvement-classic} is not feasible in this case at step $i=2$. \autoref{algo:PI} can thus not proceed although $f$, $\ell$, $h^0$ and $V^0$ are continuous and $\mathcal{U}(x)$ is compact and independent of $x$.

To overcome this issue, we present in the next section a modification of PI that ensures recursive feasibility for the above example.

\begin{rem}
The conditions in   \cite{BertsekasVIPI} for the feasibility of \eqref{eq:PI-ith-improvement-classic} are not satisfied in this example. Indeed, \cite{BertsekasVIPI} requires the set $\mathcal{U}^2_\lambda(x)=\{u \in \mathcal{U}(x)\,:\, \ell(x, u) + V^{1}(f(x,u))\leq \lambda \}$ to be compact  for any $x,\lambda\in\R$, see the discussion following  \cite[(7) and Proposition 3]{BertsekasVIPI}. However, $\mathcal{U}^2_\lambda(\bar x+1)$ is not compact for any $\lambda\geq\frac{681}{28}$, as seen in \autoref{fig:unfeasible}. We also note that the other condition for feasibility provided in \cite[Propositions 3-4]{BertsekasVIPI},  namely   $\{u\in\mathcal{U}(x)\,:\,\ell(x,u)\leq\lambda\}$ compact for any $x,\lambda\in\R$,  is verified for the considered example but does not guarantee feasibility here. \mbox{}\hfill$\Box$
\end{rem}

\begin{rem}
Contrary to model predictive control problems where the main obstacles for recursive feasibility are state constraints, we see via this example that the issue arises with PI even when no restriction is imposed on the set where the state lies. \mbox{}\hfill$\Box$
\end{rem}

\section{Policy iteration plus} \label{sec:mainresults} 
In this section we modify PI to enforce recursive feasibility  under mild conditions,  and we call the modified algorithm \PIplus.  We show \PIplus preserves the recursive (robust) stability and near-optimality properties  stated for PI in \autoref{sec:sys-cost-problem}.

\subsection{The algorithm} \label{subsec:piNew}
\PIplus is presented in \autoref{algo:PIplus}. At any iteration $i\in\Zo$, to enforce the existence of minimizer to $\ell(x,u) + V^{i}(f(x,u))$ over $u\in\mathcal{U}(x)$  in \eqref{eq:PI-ith-improvement-classic} for any given $x\in\R^{n_x}$, we first  regularize the set-valued map $H^{i+1}$   in \eqref{eq:PI-ith-improvement-reg}, see \cite[Def. 4.13]{Goebel-Sanfelice-Teel-book}. For $x\in\R^{n_x}$, the set $H_r^{i+1}(x)$ is the intersection for all $\delta>0$ of the closures of sets 
$H^{i+1}(x+\delta \mathbb{B})$.
As a result,  $H_{r}^{i+1}$ in \eqref{eq:PI-ith-improvement-reg} is  outer semicontinous\footnote{See \cite[Definition 5.4]{Rockafellar-Wets-book}.} \cite[Lemma 5.16]{Goebel-Sanfelice-Teel-book}. It is important to notice that $H^{i+1}(x)\subseteq H_{r}^{i+1}(x)$ for any $x\in\R^{n_x}$ and $i\in\Zo$.

\begin{algorithm}[t]
\renewcommand{\hypcapspace}{2\baselineskip}
\capstart
\renewcommand{\hypcapspace}{0.5\baselineskip}
 \caption{\label{algo:PIplus} Policy Iteration Plus (\PIplus)}
 \begin{algorithmic}[1]
 \renewcommand{\algorithmicrequire}{\textbf{Input:}}
 \renewcommand{\algorithmicensure}{\textbf{Output:}}
 \REQUIRE  $f$ in \eqref{eq:sys-plant}, $\ell$ in \eqref{eq:J}, initial policy $h^0\in\mathcal{U}$ 
 \ENSURE  Policy $h_r^{\star,\infty}$, cost $V_r^\infty$

 \STATE \textbf{Initial evaluation step:}  for all $x\in\R^{n_x}$,   
  \begin{equation}V_r^0(x) := J(x,h^0). \tag{\PIplus.1} \label{eq:V0-PIplus}
 \end{equation}
 \STATE \textbf{Let} $H^0:=H_r^0:=H_r^{\star,0}:=\{h^0\}$. \label{line:H0}
  \FOR  {$i\in\Zo$} 
  \STATE \textbf{Policy improvement step:}   for all $x\in\R^{n_x}$,   
  \begin{equation}
     H^{i+1}(x) := \argmin\limits_{u\in\mathcal{U}(x)} \left\{\ell(x,u)+V_r^{i}(f(x,u))\right\}.
 \tag{\PIplus.2} \label{eq:PI-ith-improvement}
\end{equation}
 \STATE \textbf{Policy regularization step:}  for all $x\in\R^{n_x}$, 
 \begin{equation}
         H_{r}^{i+1}(x) := \bigcap_{\delta>0} \overline{H^{i+1}(x+\delta \mathbb{B})}. \tag{\PIplus.3} \label{eq:PI-ith-improvement-reg} 
 \end{equation}
  \STATE \textbf{Policy evaluation step:}  for all $x\in\R^{n_x}$,  
 \begin{equation}
V_r^{i+1}(x):=\min_{h_{r}^{i+1}\in H_{r}^{i+1}} J(x,h_{r}^{i+1}).
\tag{\PIplus.4}
\label{eq:PI-ith-evaluation}
\end{equation}
\STATE \textbf{Select} $h_r^{\star,i+1}\in H_r^{\star,i+1}$ where, for all $x\in\R^{n_x}$,
\begin{equation}
  \hspace{-0.1em} {\small  H_r^{\star,i+1}(x) {:=}\!\!\!\argmin_{u\in H_r^{i+1}(x)} \!\!\!\!\!\!\left\{\ell(x,u) + V_r^{i+1}(f(x,u))\right\}\!.}  \tag{\PIplus.5}  \label{eq:PI-ith-improvement-best}
\end{equation} \vspace{-2ex}
 \ENDFOR 
 \STATE \textbf{return} $h_r^{\star,\infty} \in H_r^{\star,\infty}$ and $V_r^\infty$.
 \end{algorithmic}
\end{algorithm}

The second modification compared to PI is on the evaluation step in \eqref{eq:PI-ith-evaluation}. Given $i\in \Zo$, for any $x\in\R^{n_x}$, instead of \eqref{eq:PI-ith-evaluation-classic}, we define $V_r^{i+1}(x)$ as the minimum cost over \emph{all} selections $h_{r}^{i+1}$ of $H_{r}^{i+1}$, see \eqref{eq:PI-ith-evaluation}.
Note that all selections do not necessarily lead to the same cost $V^{i+1}$, see \autoref{ssec:counterexample} for an example.
The differences with the evaluation step in PI are that we consider $H_r^{i+1}$, instead of $H^{i+1}$, and   that we do not take an arbitrary selection of this set-valued map, but only those policies which give the minimum cost, see \eqref{eq:PI-ith-improvement-best}. Therefore,   $V_r^{i+1}(\cdot)=J(\cdot,h_r^{\star,i+1})$ for any $h_r^{\star,i+1}\in H_r^{\star,i+1}$.

We will show that these two modifications are essential to ensure the recursive feasibility of \PIplus. It can also be already noted that, when $H^{i}$ is non-empty,  outer semicontinuous and    single-valued   at any iteration $i\in\Zo$, \PIplus reduces to PI; we elaborate more on the links between PI and \PIplus in \autoref{sssec:robustPI}.
\begin{rem}
Regarding the issue identified in the example of \autoref{ssec:counterexample}, \PIplus produces cost $V^{1'}$ and selects $h^{1'}$ due \eqref{eq:PI-ith-evaluation} and \eqref{eq:PI-ith-improvement-best}, which ensures  $u\mapsto\ell(\bar{x}+1,u)+V_r^{1}(f(\bar{x}+1,u))$ is lower semicontinuous for $u\in [-\delta,1]$ and thus $H^2(\bar{x}+1)$ is non-empty. \mbox{}\hfill$\Box$
\end{rem}

\begin{rem}
In this work, we do not consider the possible errors arising from solving \eqref{eq:PI-ith-improvement-reg}, \eqref{eq:PI-ith-evaluation} and \eqref{eq:PI-ith-improvement-best}  in practice. These relevant questions are left for future work. \mbox{}\hfill$\Box$
\end{rem}

\begin{rem}
Lower semicontinuous regularizations of the optimal control problem is often considered in the literature for the optimal control of continuous-time systems, see, e.g.,  \cite{fleming2006controlled}. If we would consider a lower semicontinuous regularization of  $(x,u)\mapsto\ell(x,u)+V^i(f(x,u))$ in \eqref{eq:PI-ith-improvement-classic}, instead of \eqref{eq:PI-ith-improvement-reg}, to   guarantee $H^{i+1}(x)$ is non-empty for any $x\in\R^{n_x}$, it is then unclear here whether: (i) $V^\star\leq V^{i+1}\leq V^{i}$; and (ii) that  $V^{i+1}$ is the induced cost of some policy $h^{i+1}$. Both are key properties  that allow for $V^{i}$ to converge to $V^\star$ when $i\to\infty$. \mbox{}\hfill$\Box$
\end{rem}

\subsection{Desired properties} \label{subsec:piNew-objectives}

Like in \autoref{ss:desiredproperties}, the objectives in this section are to provide conditions under which \PIplus is such that at any iteration $i\in\Zo$:

\begin{itemize}
    \item (\emph{recursive feasibility}) for any   $x\in\R^{n_x}$,  $H^{i}(x)$, $H_{r}^{i}(x)$ and $H_{r}^{\star,i}(x)$ are non-empty;
    \item (\emph{recursive stability}) system~\eqref{eq:sys-plant} whose inputs are generated by \PIplus, i.e.,
\begin{equation}
      x(k+1) \in f(x(k),H_{r}^{\star,i}(x(k))) =: F_{r}^{i}(x(k)),
\label{eq:ith-auto-sys}
\end{equation} exhibits desirable set stability properties;
    \item (\emph{near-optimality guarantees}) explicit bounds on  $V_{r}^{i}(x)-V^{\star}(x)$ for any $x\in\R^{n_x}$ can be derived, which  asymptotically go to zero as $i$ increases.
\end{itemize}

For convenience, we abuse notation to denote solutions to system \eqref{eq:ith-auto-sys}    in the sequel as    $\phi^{i}(\cdot,x)$ when initialized at some $x\in\R^{n_x}$ for any $i\in\Zo$. As for PI, the robustness of the stability properties is deferred to \autoref{sec:robust}.

\subsection{Recursive feasibility} \label{subsec:rfeasability}
The  recursive feasibility analysis rests on the next assumptions in addition to the ones from \autoref{subsec:assumptions}, which are all checkable a priori.

\begin{ass}
The following holds.
\begin{enumerate}
\item[(i)] The function $f$ and the stage cost function $\ell$ are continuous on $\R^{n_x}\times\R^{n_u}$.
\item[(ii)] $\ell$ is level-bounded in $u$, locally uniform in $x$.
\item[(iii)] $\mathcal{U}$ is outer semicontinuous on $\R^{n_x}$.
\item[(iv)] $V_r^{0}=J(x,h^0)$ is lower semicontinuous on $\R^{n_x}$.
\item[(v)] For all $\Delta>0$, the set $\lbrace x\in\R^{n_x}:\sigma(x)\leq\Delta\rbrace$ is compact.
\mbox{} \hfill$\Box$
\end{enumerate}
\label{ass:feasibility}
\end{ass}

Item (i) of  \autoref{ass:feasibility} imposes regularity conditions on 
model $f$ and  stage cost $\ell$. Item (ii) of \autoref{ass:feasibility} is 
satisfied when $\ell(x,u)=\ell_1(x,u)+\ell_2(u)$ for any $(x,u)\in\R^{n_x}\times\R^{n_u}$ with $\ell_1(x,u)\geq0$ and $\ell_2$ radially unbounded, i.e., $\ell_2(u)\to+\infty$ as $|u|\to+\infty$ for instance. A typical example of $\ell_2$  radially unbounded being $\ell_2(u)=u^\top R u$ with $R\in\R^{n_u\times n_u}$ symmetric and positive definite; note that this property also trivially holds when $\mathcal{U}$ is independent of $x$ and compact. On the other hand, item (iii) of \autoref{ass:feasibility} is satisfied when $\mathcal{U}(\cdot)=\R^{n_u}$ or $\mathcal{U}(\cdot)$ is independent of $x$ and compact for example.  Item (iv)  of \autoref{ass:feasibility}  is a mild regularity assumption on the initial cost function, which holds when $V_r^0$ is continuous on $\R^{n_x}$. Finally, item (v)  of \autoref{ass:feasibility} implies $\sigma$ is radially unbounded, which implies that the attractor $\{x\in\R^{n_x}\,:\,\sigma(x)=0\}$ is compact.

\begin{rem}
The example   in \autoref{ssec:counterexample} verifies  \autoref{SA:well-posed}\hyperref[SA:stabilizing-initial-policy]{--\ref*{SA:stabilizing-initial-policy}} and \autoref{ass:feasibility}  with $\sigma(\cdot)=|\cdot|$, $W=\overline{\alpha}_W=0$,  $\alpha_W=\chi_W=\mathbb{I}$ and $\overline{\alpha}_V=V_r^0 \circ \sigma$ where $V_r^0=V^0$. As a consequence, the objectives stated for \PIplus in \autoref{subsec:piNew-objectives} are satisfied  for this example in view of the results presented next.
\mbox{}\hfill $\Box$
\end{rem}

The next theorem ensures  recursive feasibility for \PIplus.  For technical reasons, we prove recursive feasibility in conjunction with recursive stability for \PIplus in \refappendix{appendix:recursive}.

\begin{thm} \label{thm:rfeasability} Suppose \autoref{ass:feasibility} holds. Sets $H^i(x)$, $H_r^i(x)$ and $H_r^{\star,i}(x)$ are non-empty for any $i\in\Zo$ and $x\in\R^{n_x}$. \mbox{}\hfill$\Box$
\end{thm}

Given the recursive feasibility of \PIplus in \autoref{thm:rfeasability}, we now establish that the properties of PI in \autoref{ssec:piprop} also hold for \PIplus.

\subsection{Recursive stability}\label{subsec:rstability}

The next theorem ensures recursive stability for \PIplus. Its proof is given in \refappendix{appendix:recursive}.  
\begin{thm}\label{thm:rstability}
Suppose \autoref{ass:feasibility} holds. For any $i\in\Zo$, for any $x\in\R^{n_x}$, solution $\phi^i$ to \eqref{eq:ith-auto-sys} and $k\in\Zo$, \begin{equation}\label{eq:PIplusstab}
        \sigma(\phi^{i}(k,x))\leq\beta(\sigma(x),k)
    \end{equation} 
    where\footnote{See \autoref{foot:non-decreasing} in \autopageref{foot:non-decreasing} in \autoref{ssec:piprop}.} $\beta:(k,s)\mapsto\underline{\alpha}_Y^{-1} \circ (\mathbb{I}-\widetilde{\alpha}_Y)^{(k)}\circ \overline{\alpha}_Y(s)\in\KL$ with $\underline{\alpha}_Y,\overline{\alpha}_Y,\widetilde{\alpha}_Y$ in Table \ref{table:functions-theorems}.   
  \mbox{}\hfill $\Box$
\end{thm}

Similar to \autoref{thm:rstability-PI}, \autoref{thm:rstability} is a $\KL$-stability property of system~\eqref{eq:ith-auto-sys} with respect to $\sigma$, with $\beta$ in \eqref{eq:PIplusstab}. As previously mentioned, when $\sigma$ is a proper indicator of a compact set $\mathcal{A}$, \eqref{eq:PIplusstab} ensures that this set is uniformly globally asymptotically stable.

As in  \autoref{cor:expresult-PI}, when some of the functions in \autoref{SA:well-posed}--\ref{SA:stabilizing-initial-policy} satisfy stronger conditions, the stability property in \autoref{thm:rstability} becomes exponential. A sketch of the proof is given in \autoref{appendix:exponential}.

\begin{cor} \label{cor:expresult}
Suppose \autoref{ass:feasibility} holds and that there exist  $c_W,a_W,\overline a_V>0$ and $\overline a_W\geq0$ such that $\chi_W(s)\leq c_W s$, $\alpha_W(s)\geq a_W s$, $\overline\alpha_V(s)\leq \overline a_V s$, $\overline\alpha_W(s)\leq \overline a_W s$ for any $s\geq 0$, where $\chi_W,\alpha_W,\overline\alpha_W$ come from \autoref{SA:detect-control} and $\overline\alpha_V$ comes from \autoref{SA:stabilizing-initial-policy}. Then,
\autoref{thm:rstability} holds with $\beta:(s,k)\mapsto\dst\tfrac{\overline a_Y}{\underline a_Y}(1-\widetilde{a}_Y)^k s\in{\exp}{-}\KL$, and $\widetilde{a}_Y,\overline a_Y,\underline a_Y$ in Table \ref{table:functions-theorems}.
\mbox{} \hfill $\Box$
\end{cor}

\subsection{Near-optimality properties}\label{subsec:near-opti}

Similar to \autoref{thm:near-opti-PI}, the next theorem provides near-optimality guarantees for \PIplus. Its proof is given in the \autoref{appendix:near-opti}.  

\begin{thm} \label{thm:near-opti} Suppose \autoref{ass:feasibility} holds. For any $i\in\Zo$ and $x\in\R^{n_x}$, 
\begin{equation} \label{eq:abstract-near}
\begin{array}{rll}
    (V_r^i-V^\star)(x) & \leq & (V_r^0-V^\star)(\phi(i,x,h^\star)),
\end{array}\end{equation}
     where  $h^\star \in H^\star$ from \eqref{eq:Hstar}. Moreover, for any $i\in\Zo$ and $x\in\R^{n_x}$
     \begin{equation}
     \begin{array}{rll}
         V_r^i(x)-V^\star(x) \leq \widetilde{\alpha}\left(\beta(\sigma(x),i)\right), \label{eq:explicit-near}
         \end{array}
     \end{equation}
     where $\beta\in\KL$ comes from \autoref{thm:rstability} and $\widetilde{\alpha}=\max_{\hat{s}\in[0,s]} \widehat{\alpha}(\hat{s})$ is non-decreasing and positive definite function with  $\widehat{\alpha}$ in Table \ref{table:functions-theorems}. 
    \mbox{}\hfill $\Box$
\end{thm}

Similar to \autoref{thm:near-opti-PI},  \autoref{thm:near-opti} provides near-optimality guarantees to \PIplus. We conclude that $V_r^i \to V^\star$ as $i\to\infty$ point-wisely, as  $\beta\in\KL$ and $\lim_{i\to\infty}\sigma(\phi(i,x,h^\star))=0$ by following similar lines as in the proof of \autoref{thm:rstability}. This qualitative property of \PIplus is strengthened in \eqref{eq:explicit-near} by employing a conservative bound for \eqref{eq:abstract-near} in the form of $\widetilde{\alpha}\left(\beta(\sigma(x),i)\right)$, whose formulas for $\beta$ and $\widetilde{\alpha}$   are given in \hyperref[thm:rstability]{Theorems \ref*{thm:rstability}} \hyperref[thm:near-opti]{and \ref*{thm:near-opti}} in terms of functions in Table \ref{table:functions-theorems} and are thus known. Again, the fact that function $\beta$ in \eqref{eq:PIplusstab}  is independent of the number of iterations $i$ ensures that the upper-bound  converges to zero as $i$ increases to infinity. 

\begin{rem} As in  \autoref{rem:cstop},  the bound in  \eqref{eq:explicit-near} can be used to design stopping criteria for \PIplus to ensure a desired level of near-optimality.  \mbox{}\hfill$\Box$ 
\end{rem}

A direct consequence of \autoref{thm:near-opti} is that the sequence of cost functions $V_r^i$ converge \emph{compactly} to $V^\star$, as formalized next.

\begin{prop} \label{prop:uniform}
Suppose \autoref{ass:feasibility} holds. The sequence of functions $V_r^i$ monotonically compactly converges to $V^\star$, i.e.:
\begin{itemize}
    \item[(i)] for all $x\in\R^{n_x}$, $V_r^{i+1}(x)\leq V_r^{i}(x)$;
    \item[(ii)] for any  compact set $K\subset\R^{n_x}$ and $\delta>0$ there exists $i^\star\in\Zp$ such that for any $x\in K$ and $i\geq i^\star$, $V_r^i(x)-V^\star(x) \leq  \delta$ holds. 
$\mbox{}\hfill\Box$
\end{itemize}
\end{prop}

\autoref{prop:uniform} establishes the uniform convergence of $V_r^i$ to $V^\star$ on any compact set $K\subset\R^{n_x}$, which, compared to \autoref{prop:uniform-PI}, comes from the fact that level sets of $\sigma$ are compact in view of item (iv) of \autoref{ass:feasibility}.

\section{Robust stability} \label{sec:robust}

So far, we have established stability properties for system \eqref{eq:sys-plant}, whose inputs are generated either by PI or \PIplus. However, it is not clear yet whether these stability properties are robust to uncertainties and perturbations. Indeed, the corresponding closed-loop systems are given by difference inclusions and it has been shown in \cite{Grimm-et-al-aut04(examples)} in the context of model predictive control that difference inclusion may have zero robustness despite them satisfying so-called $\KL$-stability properties like in \hyperref[thm:rstability-PI]{Theorems \ref*{thm:rstability-PI}} \hyperref[thm:rstability-PI]{and \ref*{thm:rstability}}, in the sense that a vanishing arbitrarily small perturbation, may destroy asymptotic stability properties established for the unperturbed system. It is thus vital that the stability properties stated in \autoref{subsec:rstability} come with some nominal robustness properties. This is the focus of this section. The results are derived for system \eqref{eq:ith-auto-sys} controlled by \PIplus, but the same results follow when the inputs are generated by PI as explained in \autoref{sssec:robustPI}.

\subsection{Nominal robustness definition} \label{sssec:robustdefinition}

We consider the notion of nominal robustness in \cite{Kellett-Teel-siam-jco-05}, which we recall below in the context of this paper.
Let $\mathcal{X}\subseteq\R^{n_x}$ be open and such that  $F_{r}^i(x) \subseteq \mathcal{X}$, where $F_r^i$ comes from \eqref{eq:ith-auto-sys}, for any $x\in\mathcal{X}$ and $i\in\Zo$. We introduce a continuous function $\rho: \mathcal{X} \to \Rlo$  to perturb the set-valued map $F_r^i$ in \eqref{eq:ith-auto-sys} for $i\in\Zo$ as 
\begin{equation} \label{eq:rho-perturbed}
    \begin{split}
           \MoveEqLeft F_{r,\rho}^i(x):= \\
   &\!\!\!\!\!\left\{ \eta\in\R^{n_x}:\eta\in\{\upsilon\}+\rho(\upsilon)\mathbb{B}, \upsilon\in F_r^i(x+\rho(x)\mathbb{B})\!\right\}\!.
    \end{split}
\end{equation}
We say that $F_{r,\rho}^i$ is the $\rho$-perturbation of $F_r^i$. Given $F_{r,\rho}^i$, we have the next difference inclusion
\begin{equation} \label{eq:ith-auto-sys-perturbed}
x(k+1) \in  F_{r,\rho}^i(x(k)).
\end{equation}
System \eqref{eq:ith-auto-sys-perturbed} corresponds to system \eqref{eq:ith-auto-sys} perturbed by $\rho$, in the sense of \eqref{eq:rho-perturbed}. Note that  $F_{r}^i=F_{r,\bm{0}}^i$ and $F_{r}^i\subseteq F_{r,\rho}^i$. We denote solutions to system \eqref{eq:ith-auto-sys-perturbed} as $\phi^i_\rho(\cdot,x)$ when initialized at  $x\in\R^{n_x}$ for  $i\in\Zo$.

\autoref{thm:rstability} establishes a $\KL{-}\text{stability}$ property with respect to $\sigma$ for which the attractor is $\mathcal{A}=\{x\in\R^{n_x}\,:\,\sigma(x) = 0 \}$ when $\sigma$ is a proper indicator of $\mathcal{A}$. To define robust stability as in \cite{Kellett-Teel-siam-jco-05}, we  consider instead a two-measure $\KL{-}\text{stability}$ property with respect to $(\sigma_1,\sigma_2)$ where  $\sigma_1,\sigma_2:\mathcal{X}\to\Rlo$ are continuous. Given $\rho$ in \eqref{eq:rho-perturbed}, the attractor of perturbed system \eqref{eq:ith-auto-sys-perturbed} becomes
 \begin{equation}\label{eq:attractor-perturbed}
      \mathcal{A}^i_\rho:=\big\{\xi\in\mathcal{X} : \sup_{k\in\Zo} \sigma_1(\phi^i_\rho(k,\xi))=0 \big\}.
  \end{equation}

We are ready to define robust \KL-stability.
\begin{defn}[\cite{Kellett-Teel-siam-jco-05}]\label{def:klrobust}
Given $i\in\Zo$ and $\sigma_1,\sigma_2:\mathcal{X}\to\Rlo$ continuous, system \eqref{eq:ith-auto-sys} is  robustly \KL-stable with respect to $(\sigma_1,\sigma_2)$ on  $\mathcal{X}$   if
there exists a continuous function $\rho:\mathcal{X}\to\Rlo$ such that the following holds.
\begin{itemize}
    \item[(i)] For all $x\in\mathcal{X}$, $\{x\}+\rho(x)\mathbb{B}\subset\mathcal{X}$.
    \item[(ii)] For all  $x\in\mathcal{X}\setminus\mathcal{A}^i_{\bm{0}}$, $\rho(x)>0$.
    \item [(iii)] $\mathcal{A}^i_\rho=\mathcal{A}^i_{\bm{0}}$.
    \item [(iv)] System \eqref{eq:ith-auto-sys-perturbed} is \KL-stable with respect to $(\sigma_1,\sigma_2)$ on $\mathcal{X}$, i.e., there exists $\beta\in\KL$ such that $\sigma_1(\phi^i_\rho(k,x))\leq \beta(\sigma_2(x),k)$
     for every $x\in\mathcal{X}$ and $k\in\Zo$. When $\sigma_1\,{=}\,\sigma_2\,{=}\,\sigma$,  \eqref{eq:ith-auto-sys} is $\KL{-}\text{stable}$ with respect to $\sigma$. \mbox{}\hfill$\Box$ 
\end{itemize}     
\end{defn}
Item (i) of \autoref{def:klrobust}  imposes a condition on $\rho$ so that $F_{r,\rho}^i(\mathcal{X})\subseteq\mathcal{X}$ for any $i\in\Zo$, as $F_r^i(x)\subseteq\mathcal{X}$ for all $x\in\mathcal{X}$. Item (ii) of \autoref{def:klrobust}   in turn requires that $\rho$ is non-zero outside the attractor $\mathcal{A}^i_{\bm{0}}$ corresponding to \eqref{eq:attractor-perturbed} with $\rho=0$. Then, item (iii) of \autoref{def:klrobust} states that the perturbed attractor $\mathcal{A}^i_{\rho}$ is the same as the unperturbed one $\mathcal{A}^i_{\bm{0}}$. Finally, item (iv) of \autoref{def:klrobust}  says that the perturbed system is $\KL{-}\text{stable}$ with respect to $(\sigma_1,\sigma_2)$.

As advocated in \cite{Grimm-et-al-aut04(examples)}, if a difference inclusion is $\KL{-}\text{stable}$ with respect to some measures $(\sigma_1,\sigma_2)$, it is not necessarily robustly $\KL{-}\text{stable}$. In other words, the (unperturbed) system may have zero robustness. In the context of this paper, sufficient conditions to conclude robust $\KL{-}\text{stability}$
with respect to some pair $(\sigma_1,\sigma_2)$ on  $\mathcal{X}$ for system \eqref{eq:ith-auto-sys} with $i\in\Zo$ are that: 
\begin{itemize}
    \item[(a)] $F_r^i(x)$ is non-empty and compact for any  $x\in\mathcal{X}$;
    \item[(b)] the Lyapunov function used to establish the stability of the unperturbed system is continuous.
\end{itemize}

The next lemma states that item (a) above holds.
\begin{lem}\label{lem:compactness}
Suppose \autoref{ass:feasibility} holds. For any $i\in\Zo$,  $F_r^i(x)$ is compact and non-empty for any $x\in\R^{n_x}$. \mbox{}\hfill$\Box$
\end{lem}
\noindent\textbf{Proof:} Let $i\in\Zo$ and $x\in\R^{n_x}$, since $f$ is continuous by item (i) of \autoref{ass:feasibility}, the compactness of $F_r^i(x)= f\left(x,H_r^{\star,i}(x)\right)$ follows from the compactness of $H_r^{\star,i}$, which we now show. To this end, we prove that the conditions of \cite[Theorem 1.17(a)]{Rockafellar-Wets-book} are verified  by invoking similar arguments  to those employed in the proof of \autoref{lem:PI-ith-feasible} in \autoref{appendix:PI-ith-feasible}. Consider $g_r^{i}:(x,u)\mapsto \ell(x,u)+V_r^{i}(f(x,u))+\delta_{H_r^{i}(x)}(u)$. On the one hand, $g_r^{i}$ is proper and level-bounded in $u$, locally uniform in $x$, for the same reasons as stated  in the proof of \autoref{lem:PI-ith-feasible}.  On the other hand, $(x,u)\mapsto \ell(x,u)+V_r^{i}(f(x,u))$  is lower semicontinuous on $\R^{n_x}\times\R^{n_u}$, in view of the lower semicontinuity of $V_r^{i}$ on $\R^{n_x}$  from \autoref{prop:semilowercontinuity} in \refappendix{appendix:recursive} and the continuity of $f$ and $\ell$ by \autoref{ass:feasibility}. Moreover,  $\delta_{H_r^{i}(\cdot)}(\cdot)$ is lower semicontinuous on $\R^{n_x}\times\R^{n_u}$, as this holds from  \autoref{appendix:deltaosclsc} in \autoref{appendix:technical} since $H_r^{i}$ is an outer semicontinuous set-valued in view of \cite[Lemma 5.16]{Goebel-Sanfelice-Teel-book}. Thus $g_r^{i}$ is lower semicontinuous on $\R^{n_x}\times\R^{n_u}$ and all the conditions of \cite[Theorem 1.17(a)]{Rockafellar-Wets-book} hold, we can therefore apply it to conclude $H_r^{\star,i+1}(x)$ in \eqref{eq:PI-ith-improvement} is non-empty and compact at $x\in\R^{n_x}$. Since $x$  has been selected arbitrarily, the desired result follows. $\mbox{}\hfill\blacksquare$

Regarding item (b), it turns out that the   Lyapunov function used to establish $\KL{-}\text{stability}$ in \autoref{subsec:rstability}, namely $Y^i=\rho_V(V_r^i)+\rho_W(W)$, is only guaranteed to be lower semicontinuous. Indeed, while $W,\rho_V$ and $\rho_W$ are continuous, $V_r^i$ is only shown to be lower semicontinuous in \autoref{prop:semilowercontinuity} given in \refappendix{appendix:recursive}. Thus we will require  a different analysis or extra assumptions as presented next.

\subsection{Robust semi-global practical stability} \label{sssec:robustpractical}

An alternative Lyapunov function to analyse the stability of \eqref{eq:ith-auto-sys} is
$Y^\star:=\rho_V(V^\star)+\rho_W(W)$, where $V^\star$ is the optimal value function defined in \eqref{eq:Vstar}, which is known to be continuous in view of \cite[Theorem 3]{postoyan_stability_discounted} under a mild additional assumption on  $\mathcal{U}$ that can be checked a priori.
\begin{ass} \label{ass:Ucont}
Set-valued map $\mathcal{U}$ is continuous and locally bounded\footnote{See \cite[Definition 5.4]{Rockafellar-Wets-book} and \cite[Definition 5.14]{Rockafellar-Wets-book}, respectively.} on $\R^{n_x}$. \mbox{}\hfill$\Box$
\end{ass}

However, function $Y^\star$ only allows to establish a semiglobal practical $\KL{-}\text{stability}$ property under \hyperref[SA:well-posed]{SA\ref*{SA:well-posed}}\hyperref[SA:stabilizing-initial-policy]{--\ref*{SA:stabilizing-initial-policy}} for $i$ is sufficiently large so that, somehow, $V_r^i$ is sufficiently close to $V^\star$ in view of \autoref{prop:uniform}.
The next theorem provides a first robustness guarantee for the stability property established in  \autoref{thm:rstability}. Its proof is given in \autoref{appendix:proofrobust-practical-stability}.
\begin{thm}\label{thm:robust-practical-stability}
For any $\delta,\Delta>0$, let $\sigma_1=\max\{\sigma-\delta,0\}$, $\sigma_2=\sigma$, $\mathcal{X}=\{x\in\R^{n_x}\,:\, \sigma(x)<\Delta\}$ and suppose \hyperref[ass:feasibility]{Assumptions \ref*{ass:feasibility}} \hyperref[ass:Ucont]{and \ref*{ass:Ucont}} hold. There exists $i^\star\in\Zo$ such that, for any $i\geq i^\star$,
system \eqref{eq:ith-auto-sys} is robustly $\KL{-}\text{stable}$  with respect to $(\sigma_1,\sigma_2)$ 
on $\mathcal{X}$.\mbox{}\hfill$\Box$

\end{thm}

\autoref{thm:robust-practical-stability} implies that the stability property of system \eqref{eq:ith-auto-sys} established in \autoref{thm:rstability} is robust in a semiglobal practical sense after sufficiently many iterations $i$ with  $i\geq i^\star$  and $i^\star$ depends on $\delta,\Delta$. It is a semiglobal property as we  consider the set of initial conditions as
$\{x\in\R^{n_x}\,:\, \sigma(x)<\Delta\}$ instead of $\R^{n_x}$, and it is  practical as $\sigma_1$ is given by $\max\{\sigma(\cdot)-\delta,0\}$ instead of $\sigma$.

By strengthening the assumptions given in \autoref{subsec:assumptions}, we can ensure  stronger  robustness guarantees as shown in the next corollaries.

\begin{cor}
\label{cor:robust-asymp-semiglobal-stability}
For any $\Delta>0$, let  $\mathcal{X}=\{x\in\R^{n_x}\,:\, \sigma(x)<\Delta\}$ and suppose \hyperref[ass:feasibility]{Assumptions \ref*{ass:feasibility}} \hyperref[ass:Ucont]{and \ref*{ass:Ucont}} hold. There exists $i^\star\in\Zo$ such that, for any $i\geq i^\star$,
system \eqref{eq:ith-auto-sys} is robustly $\KL{-}\text{stable}$  with respect to $\sigma$
on $\mathcal{X}$ when there exist $l\in\Rlp$, $c_W,a_W,\overline a_V>0$ and $\overline a_W,\widetilde{a}\geq0$ such that,  for any $s\in[0,l)$, $\chi_W(s)\leq c_W s$, $\alpha_W(s)\geq a_W s$, $\overline\alpha_V(s)\leq \overline a_V s$, $\overline\alpha_W(s)\leq \overline a_W s$, $\widetilde{\alpha}(s)\leq \widetilde{a}s$. \mbox{}\hfill$\Box$
\end{cor}

\autoref{cor:robust-asymp-semiglobal-stability} extends \autoref{thm:robust-practical-stability} to allow $\delta=0$, thus enabling   robust semiglobal \emph{asymptotic} stability guarantees. The proof follows by the same manipulations  as in the proofs of \cite[Corollary 1]{granzotto2020finite} or \cite[Corollary 2]{Grimm-et-al-tac2005} and is thus omitted. 
The next corollary ensures a robust \emph{global} practical stability property under different conditions.
\begin{cor}
\label{cor:robust-practical-global-stability}
For any $\delta>0$, let $\sigma_1=\max\{\sigma(\cdot)-\delta,0\}$ and suppose \hyperref[ass:feasibility]{Assumptions \ref*{ass:feasibility}} \hyperref[ass:Ucont]{and \ref*{ass:Ucont}} hold. There exists $i^\star\in\Zo$ such that, for any $i\geq i^\star$,
system \eqref{eq:ith-auto-sys} is robustly $\KL{-}\text{stable}$  with respect to $(\sigma_1,\sigma)$ 
on $\mathcal{X}=\R^{n_x}$ when there exist $L\in\Rlo$, $c_W,a_W,\overline a_V>0$ and $\overline a_W,\widetilde{a}\geq0$ such that,  for any $s\in[L,\infty)$, $\chi_W(s)\leq c_W s$, $\alpha_W(s)\geq a_W s$, $\overline\alpha_V(s)\leq \overline a_V s$, $\overline\alpha_W(s)\leq \overline a_W s$,  $\widetilde{\alpha}(s)\leq \widetilde{a}s$.\mbox{}\hfill$\Box$
\end{cor}

\autoref{cor:robust-practical-global-stability} extends \autoref{thm:robust-practical-stability} to allow $\Delta=\infty$, thus indeed enabling for robust global practical stability guarantees. The proof also follows by similar manipulations   as in the proofs of \cite[Corollary 1]{granzotto2020finite} or \cite[Corollary 2]{Grimm-et-al-tac2005} and is thus omitted.
Finally, we can combine the conditions of \hyperref[cor:robust-asymp-semiglobal-stability]{Corollaries \ref*{cor:robust-asymp-semiglobal-stability}} \hyperref[cor:robust-practical-global-stability]{and \ref*{cor:robust-practical-global-stability}} to provide robust global and asymptotic stability properties. Its proof is omitted as it follows from \autoref{thm:robust-practical-stability}, \hyperref[cor:robust-asymp-semiglobal-stability]{Corollaries \ref*{cor:robust-asymp-semiglobal-stability}} \hyperref[cor:robust-practical-global-stability]{and \ref*{cor:robust-practical-global-stability}}.

\begin{cor}
\label{cor:robust-asymptotic-global-stability}
Suppose \hyperref[ass:feasibility]{Assumptions \ref*{ass:feasibility}} \hyperref[ass:Ucont]{and \ref*{ass:Ucont}}  hold. There exists $i^\star\in\Zo$ such that, for any $i\geq i^\star$,
system \eqref{eq:ith-auto-sys} is robustly $\KL{-}\text{stable}$  with respect to $\sigma$ 
on $\mathcal{X}=\R^{n_x}$ when both  conditions at the end of \hyperref[cor:robust-asymp-semiglobal-stability]{Corollaries \ref*{cor:robust-asymp-semiglobal-stability}} \hyperref[cor:robust-practical-global-stability]{and \ref*{cor:robust-practical-global-stability}} are satisfied. \mbox{}\hfill$\Box$
\end{cor}

The results of this section establishes robust stability properties only when sufficiently many iterations of \PIplus have been performed. When $V_r^0$ is sufficiently close to $V^\star$, the required number of iterations for robust stability may be as low as 0. In general, however, we might not have robust stability at the first iterations. In order to establish robustness of the global asymptotic stability of \PIplus at any iteration, we  require extra assumptions, which are presented in the next section.

\subsection{Robust global asymptotic stability} \label{sssec:robustcoincide}

It is possible to ensure robust stability at any iteration by adding an assumption on the problem, thanks to which $V_r^i$ can be shown to be continuous so that the desired robustness property follows from \cite{Kellett-Teel-siam-jco-05}. 

\begin{ass}\label{ass:samecost}
For any iteration $i\in\Zo$, for any two selections $h,h' \in H^i$, $J(\cdot, h) = J(\cdot,h')$.
\mbox{}\hfill$\Box$

\end{ass} 

\autoref{ass:samecost} implies that all selections at any given iteration have the same induced cost.  This arises, for example, when $H^i$ is single-valued at any $i\in\Zo$, which for instance occurs when the plant dynamics is linear and the cost is quadratic under mild conditions.  While verifying \autoref{ass:samecost} a priori is not trivial, we can check this condition at each iteration for
certain classes of optimal control problems, e.g., when the set of different selections of $H^i$ is finite. 

We establish the next key result that follows from  \autoref{ass:samecost}, whose proof is given in \autoref{appendix:continuous}.

\begin{prop} \label{prop:continuous}
Suppose \hyperref[ass:feasibility]{Assumptions \ref*{ass:feasibility}} \hyperref[ass:samecost]{and \ref*{ass:samecost}} hold  and $h^0\in\mathcal{U}$ be such that $J(\cdot,h^0)$ is continuous. Then, for any $i\in\Zo$, $V_r^i$ is continuous on $\R^{n_x}$. \mbox{}\hfill$\Box$
\end{prop}

\autoref{prop:continuous} gives a condition for  $V_r^i$ to be continuous. 
\begin{rem}
While \autoref{prop:continuous} is vital for the subsequent robustness analysis, its consequences are interesting  in their own right. Indeed, knowing that $V_r^i$ is continuous is useful when $f$ is unknown, which is not the case in this paper, as it allows using a variety of techniques to learn the value function and an associated policy at each iteration on any given compact set. \mbox{}\hfill$\Box$
\end{rem}
We provide the next robustness analysis for \PIplus  under \autoref{ass:samecost}. The proof is given in \autoref{appendix:proof-robust-asymptotic-stability}.

\begin{thm}\label{thm:robust-asymptotic-stability}
Suppose \hyperref[ass:feasibility]{Assumptions \ref*{ass:feasibility}} \hyperref[ass:samecost]{and \ref*{ass:samecost}} hold  and $h^0\in\mathcal{U}$ be such that $J(\cdot,h^0)$ is continuous. Then, for any $i\in\Zo$ system  \eqref{eq:ith-auto-sys} is robustly $\KL{-}\text{stable}$  with respect to measure $\sigma$  on $\mathcal{X}=\R^{n_x}$ for any $i\in\Zo$.\mbox{}\hfill$\Box$
\end{thm}

In contrast to \autoref{thm:robust-practical-stability}, \autoref{thm:robust-asymptotic-stability} shows that the stability guarantee given by  \autoref{thm:rstability} is indeed robust  at any iteration.

\subsection{Robustness for PI} \label{sssec:robustPI}

The robustness properties established for \PIplus also apply to PI. Indeed, when $F^{i}$ is compact and non-empty on $\R^{n_x}$ for every $i\in\Zo$, we can draw corresponding robustness guarantees for PI. In turn, $F^{i}$  compact and non-empty on $\R^{n_x}$ for every $i\in\Zo$ holds when \autoref{ass:feasibility} holds and $V^i$ is lower semicontinuous for every $i\in\Zo$.

Moreover, 
\PIplus reduces to PI under \autoref{ass:samecost} is satisfied, i.e., $V_r^i=V^i$  when $V_r^0=V^0$ and continuous. To see this, note that $V_r^0$ continuous implies that  $H^1$ is  outer semicontinuous, hence $H^1=H_r^1$. Then, \autoref{ass:samecost} guarantees that any selection of $H^1_r$  produces the same induced cost, hence $H^1_r=H^{\star,1}_r$. Thus $H^{\star,1}_r=H^1$ and $V_r^1=V^1$, which are also continuous in view of \autoref{prop:continuous} for $i=1$. The above reasoning allows to obtain $V_r^i=V^i$ holds for any $i\in\Zo$ by induction. In conclusion, we also endow  classical PI with robustness guarantees as in \autoref{sssec:robustcoincide} when  \autoref{ass:samecost} holds.

\vspace{-1ex}

%%%%%%%%%%%%%%%%%%%%%%%%%%%%%%%%%%%%%%%%%%%%%%%%%%%%%%%%%%%%
\section{Conclusion} \label{sec:conclusion}

We presented  conditions to ensure recursive robust set stability properties for  deterministic nonlinear discrete-time systems  whose inputs are generated by PI. We also gave novel near-optimality properties, which do not rely on a discount factor  contrary to the related works of the literature, see, e.g., \cite{Bertsekas-book12(adp),bertsekasbook2019rl}. Because PI may fail to be recursively feasible, we have then modified it to address this issue, which leads to the algorithm called \PIplus. \PIplus was shown to be recursively feasible under mild conditions and to preserve the robust stability and near-optimality properties of PI when the attractor is compact.  It will be interesting in future work to study the conservatism of the given near-optimality bounds and to extend the current robustness properties to address more general non-vanishing perturbations, possibly in a data-driven case, as recently done in \cite{Pang2022} for linear systems.

\appendices

\vspace{-1ex}

\section{Proofs of \texorpdfstring{\hyperref[thm:rfeasability]{Theorems \ref*{thm:rfeasability}}}{Theorems \ref*{thm:rfeasability}} and \ref{thm:rstability}}

\label{appendix:recursive}

\subsection{Proof outline}
We prove \hyperref[thm:rfeasability]{Theorems \ref*{thm:rfeasability}} and \ref{thm:rstability} together by invoking the next proposition.
\begin{prop}\label{prop:induction-statement}
For all $i\in\Zo$, the following holds.
\begin{itemize}
    \item[(i)] For any $x\in\R^{n_x}$, $H^i(x)$, $H_r^i(x)$ and $H_r^{\star,i}(x)$ are non-empty.
    \item[(ii)] For any $x\in\R^{n_x}$ and $k\in\Zo$, $\sigma(\phi^i(k,x))\leq\beta(\sigma(x),k)$ holds with $\beta\in\KL$ in \autoref{thm:rstability}.
    \item[(iii)] For any $h_r^{\star,i}\in H_r^{\star,i}$ and $x\in\R^{n_x}$, $V_r^i(x)=J(x,h_r^{\star,i}(x))\leq\overline{\alpha}_V(\sigma(x))$.
    \item[(iv)] $V_r^i$ is lower semicontinuous on $\R^{n_x}$. \mbox{}\hfill $\Box$
\end{itemize}
\end{prop}

Item (i) of \autoref{prop:induction-statement} corresponds to \autoref{thm:rfeasability}, item (ii) of \autoref{prop:induction-statement} in turn corresponds to \autoref{thm:rstability}, while items (iii) and (iv)
of \autoref{prop:induction-statement} are technical properties used in the proof of \autoref{prop:induction-statement}. 

\autoref{prop:induction-statement} is proved by induction in this appendix. We first analyse the base case, i.e., $i=0$. We establish  items (i)-(iv) \autoref{prop:induction-statement} for $i=0$ by virtue of  \autoref{SA:detect-control}\hyperref[SA:stabilizing-initial-policy]{--\ref*{SA:stabilizing-initial-policy}} and \autoref{ass:feasibility}, \eqref{eq:V0-PIplus} and \autoref{line:H0} of \autoref{algo:PIplus}. In particular, for $i=0$,  items (i), (iii) and (iv) of \autoref{prop:induction-statement} follow immediately from \autoref{SA:stabilizing-initial-policy}, \autoref{ass:feasibility} and  \autoref{algo:PIplus}, while  item (ii) of \autoref{prop:induction-statement}  is  a consequence of \autoref{SA:detect-control} and \autoref{SA:stabilizing-initial-policy}. Afterwards, we will proceed with the induction step, namely we consider $i\in\Zo$ and we assume that (i)-(iv) of \autoref{prop:induction-statement} hold for $i$, we then show that they also hold for $i+1$. We proceed in steps. First  feasibility of the improvement step is proved,  in the sense that $H^{i+1}(x)$ and $H^{i+1}_r(x)$ are non-empty for any $x\in\R^{n_x}$,  and we also show $H^{i+1}$ in \eqref{eq:PI-ith-improvement} and $H_r^{i+1}$  in \eqref{eq:PI-ith-improvement-reg} are locally bounded on $\R^{n_x}$ by virtue of the inductive hypothesis. Next, we study  the feasibility of the evaluation step, and show that $V_r^{i+1}$ in \eqref{eq:PI-ith-evaluation} is well defined and is upper-bounded by $V_r^i$, hence by $\overline{\alpha}_V\circ\sigma$. This in turn allows us to provide stability guarantees for system \eqref{eq:ith-auto-sys} for $i+1$, which is given by \eqref{eq:sys-plant} in closed-loop with $H_r^{\star,i+1}$ from \eqref{eq:PI-ith-improvement-best}. As a consequence, we prove that $V_r^{i+1}$ is lower semicontinuous, and conclude that items (i)-(iv) of \autoref{prop:induction-statement} hold for $i+1$  thereby completing the proof.

\subsection{Proof of \autoref{prop:induction-statement}}

We first verify the base case, which is that items (i)-(iv) of \autoref{prop:induction-statement} hold for $i=0$, and we then  address the induction step.

\subsubsection{Base case \texorpdfstring{$(i=0)$}{(i = 0)}} \label{sssec:piInitalize}

In view of \autoref{SA:stabilizing-initial-policy}, item (iv) of \autoref{ass:feasibility}, \eqref{eq:V0-PIplus} and \autoref{line:H0} of \autoref{algo:PIplus}, we have that items (i), (iii) and (iv) of \autoref{prop:induction-statement} hold for $i=0$.

We now show that item (ii) of \autoref{prop:induction-statement} holds for $i=0$. For this purpose, we establish in the next lemma the  existence of a strict Lyapunov function for system \eqref{eq:ith-auto-sys} for $i=0$, whose proof is given in \autoref{appendix:profLyapunov}.

\begin{lem}\label{lem:lyapunov-model-based-i=0}  There exist $\underline\alpha_Y,\overline\alpha_Y,\alpha_Y,\rho_V,\rho_W\in\Kinf$ such that $Y^{0}:= \rho_V(V_r^{0}) + \rho_W(W)$  satisfies 
\begin{equation}
\begin{split}
\underline\alpha_Y(\sigma(x)) \leq   Y^{0}(x)   & {} \leq   \overline\alpha_Y(\sigma(x)) \\
Y^{0}(\upsilon) - Y^{0}(x) & {} \leq   -\alpha_Y(\sigma(x))
\label{eq:thm-lyap-model-based-i=0}
\end{split}
\end{equation}
for any $x\in\R^{n_x}$ and $v\in F_r^{0}(x)$.  \mbox{} \hfill $\Box$
\end{lem}

The next lemma formalizes the uniform global asymptotic stability property for system \eqref{eq:ith-auto-sys} at iteration $i=0$ based on \autoref{lem:lyapunov-model-based-i=0}. Its proof is given in \autoref{appendix:proofKL-stability}.

\begin{lem}\label{lem:KL-stability-i=0}  For any $x\in\R^{n_x}$, the corresponding solution $\phi^{0}$ to \eqref{eq:ith-auto-sys}   for $i=0$ verifies $\sigma(\phi^{0}(k,x))\leq\beta(\sigma(x),k)$ with $\beta\in\KL$ in \autoref{thm:rstability}. \mbox{} \hfill $\Box$
\end{lem}

As item (ii) of \autoref{prop:induction-statement} holds for $i=0$ in view of \autoref{lem:KL-stability-i=0}, we have proved that \autoref{prop:induction-statement} holds at $i=0$. It is important for the following to note that $\beta$ is constructed in terms of $\overline{\alpha}_V,\overline{\alpha}_W,\chi_W$ and $\alpha_W$, so that we can establish later that $\beta$ is independent of $i$.

\subsubsection{Induction step  \texorpdfstring{$(``i\Rightarrow i+1")$}{(i => i+1)}}
We now proceed to the induction step, and assume \autoref{prop:induction-statement} holds for $i\in\Zo$. 

\newtheorem{IH}{Inductive Hypothesis}
\newcommand{\IHname}{IH}
\renewcommand{\theIH}{\unskip}
\begin{IH}[IH]\label{IH}
Items (i)-(iv) of \autoref{prop:induction-statement} hold for  $i$. \mbox{} \hfill $\Box$
\end{IH}

Under \autoref{IH}, we  show items (i)-(iv) of \autoref{prop:induction-statement} holds for $i+1$.

%%%%%%%%%%%%%%%%%
\noindent\emph{\textbf{Feasibility of the improvement step:}}\label{subsec:piInductionStep.1}
The following lemma ensures that $H^{i+1}$ in \eqref{eq:PI-ith-improvement} and $H_{r}^{i+1}$ in \eqref{eq:PI-ith-improvement-reg}  are non-empty and locally bounded on $\R^{n_x}$. 
Its proof is given in \autoref{appendix:PI-ith-feasible}.

\begin{lem} Set-valued maps $H^{i+1}$ in \eqref{eq:PI-ith-improvement} and $H_{r}^{i+1}$ in \eqref{eq:PI-ith-improvement-reg} are non-empty and locally bounded on $\R^{n_x}$. \mbox{} \hfill$\Box$\label{lem:PI-ith-feasible} 
\end{lem}  

We now seek to evaluate $V_r^{i+1}$.

\noindent\emph{\textbf{Feasibility of the evaluation step:}}
For the $(i+1)^\text{th}$ evaluation step of \PIplus in \eqref{eq:PI-ith-evaluation}, we need to show that the minimum cost over any selection of $H_r^{i+1}$ is well-defined. In other words, that $H_r^{\star,i+1}(x)$ in \eqref{eq:PI-ith-improvement-best} is non-empty for any $x\in\R^{n_x}$. For the sake of convenience, we extend the state vector as $\chi:=(x,u)\in\R^{n_x}\times\R^{n_u}$ and we write  \begin{equation}
\chi(k+1)\in \big(f(x,u),H_{r}^{i+1}\big(f(x,u)\big)\big)=:G^{i+1}(\chi(k)).\label{eq:Gi+1}
\end{equation} 
Given any $x\in\R^{n_x}$, we consider initial conditions to \eqref{eq:Gi+1} of the form $(x,u)$ with $u\in H_{r}^{i+1}(x)$ with $H_{r}^{i+1}$ defined in \eqref{eq:PI-ith-improvement-reg}, and we denote the associated set of solutions to \eqref{eq:Gi+1} by $\mathcal{S}^{i+1}(x)$, which only depends on $x$. To extend the state vector as in~\eqref{eq:Gi+1} allows  to write the stage cost as a function of $\chi$ only, i.e., $\ell(x,u)=\ell(\chi)$, which is convenient in the sequel.
We now evaluate the cost of solutions to \eqref{eq:Gi+1}, that is, in view of \eqref{eq:PI-ith-evaluation}, for any $x\in\R^{n_x}$,
\begin{equation}
V_r^{i+1}(x)=\min_{\psi^{i+1}\in\mathcal{S}^{i+1}(x)}\sum_{k=0}^\infty\ell\big(\psi^{i+1}(k,x)\big),\label{eq:Vi+1}
\end{equation}
where $\psi^{i+1}=(\phi^{i+1},v^{i+1})\in\mathcal{S}^{i+1}(x)$ is a  solution to system \eqref{eq:Gi+1} initialized at $(x,h_{r}^{i+1}(x))$ with $h_{r}^{i+1}(x)\in H_{r}^{i+1}(x)$.  To guarantee the existence of a minimum in \eqref{eq:Vi+1}, we invoke similar arguments as in \cite[Claim 24]{Kellett-Teel-scl04}. For this purpose, we first state the next properties of $G^{i+1}$ in \eqref{eq:Gi+1}.
\begin{lem}
Function $G^{i+1}$ in \eqref{eq:Gi+1} is outer semicontinuous, locally bounded on $\R^{n_x}\times\R^{n_u}$ and $G^{i+1}(\chi)$ is non-empty for any $\chi\in\mathcal{W}^{i+1}:=\{(x,u)\in\R^{n_x}\times\R^{n_u}\,: u\in H_{r}^{i+1}(x) \}$. \mbox{} \hfill$\Box$\label{prop:osc-Gi+1}
\end{lem}

\textbf{Proof:} Since single-valued  $f$  is continuous on $\R^{n_x}\times\R^{n_u}$ by item (i) of \autoref{ass:feasibility}, it is locally bounded and outer-semicontinuous on $\R^{n_x}\times\R^{n_u}$ by \cite[Corollary 5.20]{Rockafellar-Wets-book}. Moreover, as $H_{r}^{i+1}$ in~\eqref{eq:PI-ith-improvement-reg} is locally bounded by \autoref{lem:PI-ith-feasible} and outer semicontinuous on $\R^{n_x}$ by  \cite[Lemma 5.16]{Goebel-Sanfelice-Teel-book}, we deduce that $G^{i+1}$ is locally bounded and outer semicontinuous  on $\R^{n_x}\times\R^{n_u}$ by \cite[Proposition  5.52(a) and (b)]{Rockafellar-Wets-book}. Finally, $G^{i+1}(\chi)$ is non-empty for any $\chi\in\mathcal{W}^{i+1}$ as $H_r^{i+1}(x)$ is non-empty for $x\in\R^{n_x}$ by \autoref{lem:PI-ith-feasible}. \mbox{} \hfill $\blacksquare$

We verify the existence of a policy as in  \eqref{eq:PI-ith-improvement-best} for $i+1$, and we provide key properties of such policy in the next proposition, whose proof is in \autoref{appendix:proofVi}.

\begin{lem}\label{prop:Vi+1Vi} There exists $h_{r}^{\star,i+1}\in H_{r}^{\star,i+1}$ such that $V_r^{i+1}(x)=J(x,h_{r}^{\star,i+1})$ for all $x\in\R^{n_x}$. Furthermore, the following holds for any $x\in\R^{n_x}$,
\begin{enumerate}
\item[(i)]$V_r^{i+1}(x)=\ell(x,h_{r}^{\star,i+1}(x))+V_r^{i+1}\big(f(x,h_{r}^{\star,i+1}(x))\big)$,
\item[(ii)] $V_r^{i+1}(x)\leq V_r^{i}(x) \leq \overline{\alpha}_V(\sigma(x))$. \mbox{} \hfill $\Box$
\end{enumerate}
\end{lem} 

In view of \hyperref[lem:PI-ith-feasible]{Lemmas \ref*{lem:PI-ith-feasible}} and \ref{prop:Vi+1Vi}, items (i) and (iii) of \autoref{prop:induction-statement} are verified at $i+1$. We now establish stability properties  as in item (ii) of \autoref{prop:induction-statement} at $i+1$, which then allows to show the lower semicontinuity of $V_r^{i+1}$, which is the final item  (iv) of \autoref{prop:induction-statement}. 

%%%%%%%%
\noindent\emph{\textbf{Stability:}} 
We follow similar lines as in \autoref{sssec:piInitalize} to analyze the stability of system~\eqref{eq:ith-auto-sys} at iteration $i+1$.
The next result establish the existence of a strict Lyapunov function for the $(i+1)^\text{th}$-step of \PIplus, similar to the one found in \autoref{lem:lyapunov-model-based-i=0}. The proof relies on similar arguments as the proof of  \autoref{lem:lyapunov-model-based-i=0} in \autoref{appendix:profLyapunov}, and is hence omitted.

\begin{lem}\label{lem:lyapunov-model-based-i+1} For any $x\in\R^{n_x}$ and $v\in F_r^{i+1}(x)$,
\begin{equation}
\begin{split}
\underline\alpha_Y(\sigma(x)) \leq   Y^{i+1}(x)   & {} \leq   \overline\alpha_Y(\sigma(x)) \\
Y^{i+1}(\upsilon) - Y^{i+1}(x) & {} \leq   -\alpha_Y(\sigma(x))
\label{eq:thm-lyap-model-based}
\end{split}
\end{equation}
with $Y^{i+1}:= \rho_V(V_r^{i+1}) + \rho_W(W)$ and $\underline\alpha_Y, \allowbreak \overline\alpha_Y, \allowbreak \alpha_Y, \allowbreak \rho_V, \allowbreak \rho_W \in\Kinf$ as in \autoref{lem:lyapunov-model-based-i=0}.  \mbox{} \hfill $\Box$
\end{lem}

The next result follows from \autoref{lem:lyapunov-model-based-i+1}. Its proof is very similar to the proof of \autoref{lem:KL-stability-i=0} in \autoref{appendix:proofKL-stability}, it is therefore omitted.

\begin{lem}\label{lem:KL-stability-i+1}  For any $x\in\R^{n_x}$,  any solution $\phi^{i+1}$ to \eqref{eq:ith-auto-sys}   verifies $\sigma(\phi^{i+1}(k,x))\leq\beta(\sigma(x),k)$ with $\beta\in\KL$ in \autoref{thm:rstability}. \mbox{} \hfill $\Box$
\end{lem}
Since, in  \autoref{lem:KL-stability-i=0}, we have $\beta$  only constructed in terms  of $\overline{\alpha}_V,\overline{\alpha}_W,\chi_W$ and $\alpha_W$, we therefore obtain $\beta$  independent of $i+1$ in \autoref{lem:KL-stability-i+1}. This establishes item (ii) of \autoref{prop:induction-statement} at $i+1$. All that remains to prove is item (iv) of \autoref{prop:induction-statement} for $i+1$.

%%%%%%%%
\noindent\emph{\textbf{Lower semicontinuity of  \texorpdfstring{$V_r^{i+1}$}{V\_r\textasciicircum\{i+1\}}}:}
The next lemma establishes the lower semicontinuity of $V_r^{i+1}$, whose proof is given in \autoref{appendix:lowersemicontinuity} and is inspired by  \cite[Theorem 6]{Kellett-Teel-scl04}.
\begin{lem}
$V_r^{i+1}$ is lower semicontinuous on $\R^{n_x}$.\mbox{} \hfill$\Box$\label{prop:semilowercontinuity}
\end{lem}

We have obtained lower semicontinuous $V_r^{i+1}$ and items (i)-(iv) of \autoref{prop:induction-statement} holds for $i+1$. Therefore \autoref{prop:induction-statement} holds for $i+1$ and the induction proof of \autoref{prop:induction-statement} is complete. We then deduce that \hyperref[thm:rfeasability]{Theorems \ref*{thm:rfeasability}} and \ref{thm:rstability} hold.

The remaining part of \refappendix{appendix:recursive} is dedicated to the proofs of the various lemmas stated above to prove \autoref{prop:induction-statement}.

%%%%%%%%%%%%%
\subsection{Proof of \autoref{lem:lyapunov-model-based-i=0}}\label{appendix:profLyapunov}

The proof follows similar lines as the proofs of \cite[Theorem 1]{Grimm-et-al-tac2005} and \cite[Theorem 1]{postoyan_stability_discounted}. We distinguish two cases. 

\textbf{Case where $\bm{\chi_W\leq\mathbb{I}}$:} Let $\rho_V:=\rho_W:=\mathbb{I}$, thus  $Y^0=\rho_V(V_r^0)+\rho_W(W)=V_r^0+W$.
By \autoref{SA:detect-control} and \autoref{SA:stabilizing-initial-policy}, $Y^0(x)\leq\overline{\alpha}_V(\sigma(x))+\overline{\alpha}_W(\sigma(x))$ for any $x\in\R^{n_x}$, therefore $\overline{\alpha}_Y:= \overline{\alpha}_V+\overline{\alpha}_W$. On the other hand, by \autoref{SA:detect-control} and since $\chi_W\leq\mathbb{I}$, $\alpha_W(\sigma(x))\leq \ell(x,u)+W(x)$ for any $(x,u)\in\mathcal{W}$, hence $Y^0(x) = V_r^0(x)+W(x)\geq \ell(x,h^0(x))+W(x)\geq \alpha_W(\sigma(x))$,  therefore $\underline{\alpha}_Y:= \alpha_W$. 
Moreover, since $V_r^0(x)=J(x,h^0)=\ell(x,h^0(x))+J(f(x,h^0(x)),h^0)=\ell(x,h^0(x)) +V_r^0(f(x,h^0(x)))$, we have  for $\upsilon=f(x,h^0(x))$ that
\begin{equation}
   V_r^0(\upsilon)-V^0(x)=-\ell(x,h^0(x))\label{eq:Vdiff-chi-leq-I}.
\end{equation} Again by \autoref{SA:detect-control} as $\chi_W\leq\mathbb{I}$ and $u=h^0(x)$,
\begin{equation}
    W(\upsilon)-W(x)\leq-\alpha_W(\sigma(x))+\ell(x,u). \label{eq:Wdiff-chi-leq-I}
\end{equation} 
Summing \eqref{eq:Vdiff-chi-leq-I} and \eqref{eq:Wdiff-chi-leq-I}, we obtain
$Y^0(\upsilon)-Y(x)\leq-\alpha_W(\sigma(x))$, therefore $\alpha_Y:=\alpha_W$. The case is complete.

\textbf{Case where $\bm{\chi(s)>s}$ for some $\bm{s>0}$:} We define $q_V:=\dst 2\chi_W(2\mathbb{I})$ and $q_W:=\frac{1}{2}[\chi_W+(\overline\alpha_W+\mathbb{I})\circ\alpha_W^{-1}(2\chi_W)]^{-1}$, where the involved functions come from \autoref{SA:detect-control}. Note that $q_W$ is well-defined as $\chi_W+(\overline\alpha_W+\mathbb{I})\circ\alpha_W^{-1}(2\chi_W)$ is invertible, being a class-$\Kinf$ function. Moreover, both $q_V$ and $q_W$ are of class-$\Kinf$. We further define $\rho_V(s):=\int_{0}^{s}q_V(\tau)d\tau$ and $\rho_W(s):=\int_{0}^{s}q_W(\tau)d\tau$ for any $s\geq 0$. Note that $\rho_V$ and $\rho_W$ are also of class-$\Kinf$. We are ready to define the Lyapunov function used to prove stability.

Let $x\in\R^{n_x}$ and $v\in F_r^{0}(x)$. According to  \autoref{SA:detect-control} and \autoref{SA:stabilizing-initial-policy}, $Y^0(x)\leq \rho_V\circ\overline\alpha_V(\sigma(x))+\rho_W\circ\overline\alpha_W(\sigma(x))=:\overline\alpha_Y(\sigma(x))$ and $\overline\alpha_Y$ is of class-$\Kinf$ in view of the properties of the functions $\rho_V,\rho_W,\overline\alpha_V$ and $\overline\alpha_W$. On the other hand, from \eqref{eq:Vdiff-chi-leq-I}, for any $x\in\R^{n_x}$,
\begin{align}
V_r^{0}(x)\geq \ell(x,h^0(x))
\label{eq:value-i=0-lower-bound}
\end{align}
and in view of  \autoref{SA:detect-control}, $Y^0(x) \geq \rho_V(\ell(x,h^{0}(x))) + \rho_W\left(\max\left\{\alpha_W(\sigma(x)) - \chi_W(\ell(x,h^{0}(x)))), 0\right\}\right)$. When $\frac{1}{2}\alpha_W(\sigma(x)) \leq \chi_W(\ell(x,h^{0}(x)))$, it follows $Y^0(x) \geq  \rho_V\circ\chi_W^{-1}(\frac{1}{2}\alpha_W(\sigma(x)))$. When $\frac{1}{2}\alpha_W(\sigma(x)) \geq \chi_W(\ell(x,h^{0}(x))))$ then $Y^0(x) \geq  \rho_W \circ \frac{1}{2} \alpha_W( \sigma(x) )$. Hence $Y^0(x) \geq  \underline\alpha_Y(\sigma(x))$ with $\underline\alpha_Y:=\min\left\{\rho_V\circ\chi_W^{-1}(\frac{1}{2}\alpha_W),\rho_W(\frac{1}{2}\alpha_W)\right\}\in\Kinf$. We have proved the first line in \eqref{eq:thm-lyap-model-based-i=0}.

In view of \eqref{eq:Vdiff-chi-leq-I}, \eqref{eq:value-i=0-lower-bound} and  \autoref{SA:stabilizing-initial-policy}, the conditions of \cite[Lemma 3]{Grimm-et-al-tac2005} are satisfied with $h(x)=\ell(x,h^0(x))$, $\alpha_1=\overline\alpha_V$ and $\alpha_2=0$. Thus, we derive that
\begin{align} 
\MoveEqLeft \rho_{V}(V_r^{0}(v))  - \rho_{V}(V_r^{0}(x)) \nonumber   \\
 & \leq -\chi_W(\ell(x, h^{0}(x)))\ell(x,h^{0}(x))). \label{eq:proof-thm-lyap-prel-rho-V}
\end{align}
On the other hand, the conditions of \cite[Lemma 4]{Grimm-et-al-tac2005} are verified with\footnote{$\mathbb{I}$ is added to ensure $\alpha_1\in\Kinf$.} $\alpha_1=\overline\alpha_W + \mathbb{I}$, $\alpha_2=\chi_W$ and $\alpha_3=\alpha_W$, according to \autoref{SA:detect-control}. Hence, we deduce from \cite[Lemma 4]{Grimm-et-al-tac2005}
\begin{align}
\MoveEqLeft \rho_{W}(W(v))  - \rho_{W}(W(x)) & \nonumber\\
& \!\!\!\!\!\!\!\!{} \leq 2q_{W}\Big(\chi_W(\ell(x,h^{0}(x)))+\alpha_1\circ\alpha_W^{-1}(2\chi_W(\ell(x,h^{0}(x))))\Big)\nonumber\\
 &\!\!\!\!\!\!\!\!\qquad \times\chi_W(\ell(x,h^{0}(x)))  -q_{W}(\tfrac{1}{4}\alpha_W(\sigma(x)))\tfrac{1}{4}\alpha_W(\sigma(x))\nonumber\\
 & \!\!\!\!\!\!\!\!{} =  \ell(x,h^{0}(x))\chi_W(\ell(x,h^{0}(x))) - \widetilde\alpha_W(\sigma(x))
\label{eq:proof-thm-lyap-rho-W}
\end{align}
where $\widetilde\alpha_W=q_{W}(\frac{1}{4}\alpha_W)\frac{1}{4}\alpha_W\in\Kinf$. In view of \eqref{eq:proof-thm-lyap-prel-rho-V} and \eqref{eq:proof-thm-lyap-rho-W}, it follows that 
\iftoggle{onecolumnversion}{\begin{align*}
\MoveEqLeft Y^0(v)-Y^0(x)  \\ 
& {}\leq  -\chi_W(\ell(x,h^{0}(x)))\ell(x,h^{0}(x))) \\ & \hphantom{\leq} {}+\ell(x,h^{0}(x))\chi_W(\ell(x,h^{0}(x))) - \widetilde\alpha_W(\sigma(x))\\
& {}\leq   - \widetilde\alpha_W(\sigma(x)). \numberthis
\end{align*}}{$Y^0(v)-Y^0(x)  \leq  -\chi_W(\ell(x,h^{0}(x)))\allowbreak\ell(x,h^{0}(x)))+\ell(x,h^{0}(x))\chi_W(\ell(x,h^{0}(x))) - \widetilde\alpha_W(\sigma(x))\leq   - \widetilde\alpha_W(\sigma(x))$. }%
We have shown that the second line in \eqref{eq:thm-lyap-model-based-i=0} holds with $\alpha_Y=\widetilde\alpha_W\in\Kinf$. The case is complete.

%%%%%%%%%%%%%

\subsection{Proof of \autoref{lem:KL-stability-i=0}} \label{appendix:proofKL-stability}
The proof follows by application of \autoref{lem:lyapunov-model-based-i=0} and \autoref{lem:generatingbeta} in \refappendix{appendix:technical}. Indeed, we apply \autoref{lem:generatingbeta} with $\alpha=\mathbb{I}-\widetilde{\alpha}_Y$, where $\widetilde{\alpha}_Y:=\alpha_Y \circ \overline{\alpha}_Y^{-1}$. For any $x\in\R^{n_x}$, we consider the sequence $s_k(x)$, $k\in\Zo$, defined as $s_k(x)=Y^0(\phi(k,x))$. We note that $s_{k+1}(x)\leq \alpha(s_k(x))$ as required. This property follows from  \eqref{eq:thm-lyap-model-based-i=0}, which implies  $\widetilde\alpha_Y(Y^{0}(x))=\alpha_Y\circ\overline{\alpha}_Y^{-1}(Y^0(x))\leq\alpha_Y(\sigma(x))$ hence $Y^0(\phi(k+1,x))\leq Y^0(\phi(k,x))-\widetilde{\alpha}_Y(Y^0(\phi(k,x)))$ for any $x\in\R^{n_x}$ and $k\in\Zo$. Note moreover that $\mathbb{I}-\widetilde{\alpha}_Y$ is indeed continuous, zero at zero,  $\alpha(s)<s$ for all $s>0$ and non-negative. Hence \autoref{lem:generatingbeta} generates $\widetilde{\beta}_x\in\KL$  such that $Y^0(\phi^0(k,x))\leq \widetilde{\beta}_x(Y^{0}(x),k)$ for any $x\in\R^{n_x}$ and $k\in\Zo$, and since $\alpha$ is uniform in $x$, $\widetilde{\beta}_x=:\widetilde{\beta}$. Again by \eqref{eq:thm-lyap-model-based-i=0}, we obtain that $\underline\alpha_Y(\sigma(\phi^{0}(k,x))) \leq \widetilde{\beta}(\overline\alpha_Y(\sigma(x)),k)$. Hence the result of \autoref{lem:KL-stability-i=0} holds with $\beta=\underline\alpha_Y^{-1}(\widetilde{\beta}(\overline\alpha_Y(\cdot),\cdot))\in\KL$.

%%%%%%

\subsection{Proof of \autoref{lem:PI-ith-feasible}} \label{appendix:PI-ith-feasible}

We prove that the conditions of \cite[Theorem 1.17(a)]{Rockafellar-Wets-book} are verified to first show that $H^{i+1}(x)$ is non-empty for any $x\in\R^{n_x}$. Let $g^{i+1}:(x,u)\mapsto \ell(x,u)+V_r^{i}(f(x,u))+\delta_{\mathcal{U}(x)}(u)$ be defined on $\R^{n_x}\times\R^{n_u}$. The map $g^{i+1}$ is lower semicontinuous on $\R^{n_x}\times\R^{n_u}$ as: (i) $\ell$ and $f$ are continuous on $\R^{n_x}\times\R^{n_u}$ by item (i) of \autoref{ass:feasibility}; (ii) $V_r^{i}$ is lower semicontinuous on $\R^{n_x}$  by \autoref{IH}; (iii) $(x,u)\mapsto\delta_{\mathcal{U}(x)}(u)$ is lower semicontinuous on $\R^{n_x}\times\R^{n_u}$ by \autoref{appendix:deltaosclsc} in \refappendix{appendix:technical} as $\mathcal{U}$ is outer semicontinuous on $\R^{n_x}$ by item (iii) of \autoref{ass:feasibility}. 
Moreover, $g^{i+1}$ is proper as $g^{i+1}(x,u)>-\infty$ for any $(x,u)\in\R^{n_x}\times\R^{n_u}$ and $g^{i+1}(x,u)<\infty$ for any $(x,u)\in\mathcal{W}$, as $\mathcal{U}(x)$ is non-empty for any $x\in\R^{n_x}$. In addition, as $g^{i+1}(x,u)\geq \ell(x,u)$ for any $(x,u)\in\R^{n_x}\times\R^{n_u}$ and $\ell$ is level-bounded in $u$, locally uniform in $x$ by item (ii) of \autoref{ass:feasibility},  $(x,u)\mapsto g^{i+1}(x,u)$ is level-bounded in $u$, locally uniform in $x$ for any $x\in\R^{n_x}$ by \autoref{lem:levelboundeduniform} in \refappendix{appendix:technical}. All the conditions of \cite[Theorem 1.17(a)]{Rockafellar-Wets-book} thus hold, we can therefore apply this result to deduce that {$H^{i+1}(x)$} in \eqref{eq:PI-ith-improvement} is non-empty and compact for any $x\in\R^{n_x}$. Since $H^{i+1}(x)\subseteq H_r^{i+1}(x)$ for any $x\in\R^{n_x}$ by definition from \eqref{eq:PI-ith-improvement-reg}, $H_r^{i+1}(x)$ is also non-empty for any $x\in\R^{n_x}$.

On the other hand, $H_{r}^{i+1}$ is outer semicontinuous by \cite[Lemma 5.16]{Goebel-Sanfelice-Teel-book}.  Furthermore, again by \cite[Lemma 5.16]{Goebel-Sanfelice-Teel-book}, $H_{r}^{i+1}$ is locally bounded when $H^{i+1}$ also is, which we now show by applying  \cite[Theorem 7.41(a)]{Rockafellar-Wets-book}. Indeed, let $x\in\R^{n_x}$ and denote  $\overline{V}^{i+1}(x):=\min_{u\in\mathcal{U}(x)} \left\{\ell(x,u) + V_r^{i}(f(x,u))\right\}$, which is well-defined as $H^{i+1}$ is non-empty. Since $\overline{V}^{i+1}(x) \leq \ell(x,h_r^{\star,i}(x))+V_r^{i}(f(x,h_r^{\star,i}(x))) = V_r^{i}(x)$ for $h_r^{\star,i}\in H_r^{\star,i}$ (recall that $H_r^{\star,i}(x)$ is non-empty for any $x\in\R^{n_x}$ by \autoref{IH}) and $V_r^{i}(x)\leq\overline{\alpha}_V(\sigma(x))$ again  by the \autoref{IH}, we have $\overline{V}^{i+1}(x)\leq\overline{\alpha}_V(\sigma(x))$. $\overline{V}^{i+1}(x)$ is thus bounded from above at any $x\in\R^{n_x}$. To see this, let $\mathcal{V}$ be a neighborhood of $x$, and since $\sigma$ and $\overline{\alpha}_V$  are continuous by \autoref{SA:stabilizing-initial-policy}, there exists some finite $M>0$ such that $\overline{V}^{i+1}(x)\leq \sup_{x'\in\mathcal{V}}\{\overline{\alpha}_V(\sigma(x')) \}=M$. Hence $\overline{V}^{i+1}$ is bounded from above in a neighborhood of $x$. By invoking \cite[Theorem 7.41(a)]{Rockafellar-Wets-book}, it follows that $H^{i+1}$ is locally bounded in $\mathcal{V}$. Since $x$ is arbitrary, it follows that $H^{i+1}$ is locally bounded on $\R^{n_x}$ and so is $H_r^{i+1}$.

\subsection{Proof of \autoref{prop:Vi+1Vi}} \label{appendix:proofVi}

First, we show  there exists $h_r^{i+1}\in H_{r}^{i+1}$ such that $J(\cdot,h_r^{i+1})\leq V_r^{i}$. For this purpose, since $H^{i+1}\subseteq H_{r}^{i+1}$ by definition of $H_{r}^{i+1}$ in \eqref{eq:PI-ith-improvement-reg},  it is sufficient to show $J(\cdot,h^{i+1})\leq V_r^{i}$  for any $h^{i+1}\in H^{i+1}$. To do so, we proceed  like in \cite[Section 4.2]{Sutton-Barto-book2018}. Let $x\in\R^{n_x}$. We have $V_r^{i}(x)=J(x,h_r
^{\star,i})\leq\overline{\alpha}_V(\sigma(x))$ for any $h_r
^{\star,i}\in H_r
^{\star,i}$ by  \autoref{IH} and  let $h^{i+1}\in H^{i+1}$, which is possible in view of \autoref{lem:PI-ith-feasible}. It follows that \iftoggle{onecolumnversion}{
\begin{equation}
    \begin{split}
        V_r^{i}(x)& {} =   \ell(x,h_r
^{\star,i}(x))+V_r^{i}(f(x,h_r
^{\star,i}(x)))\\
& {} \geq  \ell(x,h^{i+1}(x))+V_r^{i}(f(x,h^{i+1}(x))).
    \end{split}
\end{equation}}{$V_r^{i}(x)  =   \ell(x,h_r
^{\star,i}(x))+V_r^{i}(f(x,h_r
^{\star,i}(x)))\geq  \ell(x,h^{i+1}(x))+V_r^{i}(f(x,h^{i+1}(x)))$. }%
By repeating the same reasoning on $V_r^{i}(f(x,h^{i+1}(x)))$ and denoting $v=f(x,h^{i+1}(x))$, we derive that 
\iftoggle{onecolumnversion}{
\begin{equation} 
    \begin{split}
        V_r^{i}(x)  &{}\geq  \ell(x,h^{i+1}(x))+V_r^{i}(v) \\
& {} =   \ell(x,h^{i+1}(x))+\ell\left(v,h_r
^{\star,i}(v)\right) \\ &\quad {} + V_r^{i}(f(v,h_r
^{\star,i}(v))) \\
& {} \geq   \ell(x,h^{i+1}(x))+\ell(v,h^{i+1}(v)) \\ & \quad {} + V_r^{i}(f(v,h^{i+1}(v))).
    \end{split}
    \label{eq:proof-prop-convergence-monotone-decrease-prel}
\end{equation}}{$V_r^{i}(x)\geq  \ell(x,h^{i+1}(x))+V_r^{i}(v)  =   \ell(x,h^{i+1}(x))+\ell\left(v,h_r
^{\star,i}(v)\right)  + V_r^{i}(f(v,h_r
^{\star,i}(v)))  \geq   \ell(x,h^{i+1}(x))+\ell(v,h^{i+1}(v)) + V_r^{i}(f(v,h^{i+1}(v)))$. }%
We can continue this process infinitely many times to obtain an infinite sum, which is bounded as such sum is: (i) upper-bounded by $V_r^{i}(x)$, which itself is finite in view of $V_r^{i}(x)\leq\overline{\alpha}_V(\sigma(x))$  by \autoref{IH}; (ii) lower-bounded by $0$, since it is a sum of non-negative terms. In this manner, we derive that
\begin{equation}
V_r^{i}(x) \geq J(x,h^{i+1})
\label{eq:proof-prop-convergence-monotone-decrease-generic}
\end{equation}
for any $h^{i+1}\in H^{i+1}\subseteq H_{r}^{i+1}$. Thus, we can apply the steps of \cite[Claim 24]{Kellett-Teel-scl04} and derive that there exists $h_{r}^{\star,i+1}\in H_{r}^{i+1}$ such that $V_r^{i+1}(x)=J(x,h_{r}^{\star,i+1})$ in view of the outer semicontinuity and local boundedness of $H_{r}^{i+1}$. Hence, $H_{r}^{\star,i+1}(x)$ is non-empty for all $x\in\R^{n_x}$. Furthermore, by definition of $V_r^{i+1}$ in \eqref{eq:PI-ith-evaluation}, we have that $V_r^{i+1}(x)\leq J(x,h_r^{i+1})$ for any $h_r^{i+1}\in H_{r}^{i+1}$, and by \eqref{eq:proof-prop-convergence-monotone-decrease-generic}, we conclude that $V_r^{i+1}(x) \leq J(x,h^{i+1}) \leq V_r^{i}(x)\leq\overline{\alpha}_V(\sigma(x))$. On the other hand, $J(x,h_{r}^{\star,i+1})=\ell(x,h_{r}^{\star,i+1}(x))+J(f(x,h_{r}^{\star,i+1}(x)),h_{r}^{\star,i+1})$   for $h_{r}^{\star,i+1}\in H_{r}^{\star,i+1}$
and hence items (i)-(ii) of \autoref{prop:Vi+1Vi} are verified and the proof is complete. 

%%%%%%%%%%%%

\subsection{Proof of \autoref{prop:semilowercontinuity}} \label{appendix:lowersemicontinuity}
To prove the lower semicontinuity of $V_r^{i+1}$, we show that, for any $\varepsilon>0$, for any $x\in\R^{n_x}$ and any sequence\footnote{The limits of sequences in the appendices are understood for $n\to\infty$.} $x_n\to x$, there exists $N\in\Zp$ such that for any $n\geq N$,  $V_r^{i+1}(x_n)\geq V_r^{i+1}(x)-\varepsilon$.

Since $H_r^{\star,i+1}(x)\subseteq H_r^{i+1}(x)$ for any $x\in\R^{n_x}$, by \autoref{lem:KL-stability-i+1} there exists  $\psi=(\phi,v)\in\mathcal{S}^{i+1}(x)$ with $\mathcal{S}^{i+1}$ defined after \eqref{eq:Gi+1}, such that for any $x\in\R^{n_x}$ and  $k\in\Zo$, 
\begin{equation}
\sigma(\phi(k,x))\leq\beta(\sigma(x),k)) \label{eq:beta}
\end{equation}  
where $\beta\in\mathcal{KL}$ as in \autoref{lem:KL-stability-i+1}. From now on, we designate such solutions as ``optimal'' and we follow similar lines as in the proof of \cite[Theorem 6]{Kellett-Teel-scl04}. Let  {$x\in\R^{n_x}$} and consider an arbitrary sequence $x_n\in\R^{n_x}$, $n\in\Zo$,  which converges to $x$, as well as an arbitrary sequence $u_n$, which converges to $u\in\R^{n_u}$ with $u_n\in H_{r}^{i+1}(x_n)$. We have $u\in H_{r}^{i+1}(x)$ by outer semicontinuity of $H_{r}^{i+1}$ in  \autoref{lem:PI-ith-feasible}.
Let $\varepsilon>0$ and $\overline N\in\Zo$ be sufficiently big such that
\begin{equation}
\beta(\sigma(x)+1,\overline N)\leq \varepsilon_1:=\min\{\frac{1}{2} {\alphav}^{-1}(\frac{\varepsilon}{4}),1\},\label{eq:eps2} 
\end{equation} with $\alphav$  from \autoref{SA:stabilizing-initial-policy}. {We assume without loss of generality that $n\in\Zo$ is sufficiently large such that $\sigma(x_n)\leq\sigma(x)+1$, which is possible since $\sigma$ is continuous and {the sequence $x_n$ converges to $x$} as $n$ tends to $\infty$.} {In view of \eqref{eq:beta}}, we have for all $k\in\Zo$ and any optimal solution $(\phi,\upsilon)\in \mathcal{S}^{i+1}(x_n)$, 
\begin{equation}
\begin{aligned}
\sigma\big(\phi(k,x_n)\big)&\leq\beta\big(\sigma(x_n),k\big)\leq \beta\big(\sigma(x_n),0\big)\\
&\leq \beta\big(\sigma(x)+1,0\big).
\end{aligned}\label{eq:phihi<beta}
\end{equation}
According to item (iv) of \autoref{ass:feasibility} {and the continuity of $\beta$ and $\sigma$}, \begin{equation} \label{eq:M}
    \mathcal{M}(x):=\lbrace z\in\R^{n_x}:\sigma(z)\leq\beta\big(\sigma(x)+1,0\big)+1\rbrace
\end{equation} is compact. {Moreover, $H_{r}^{i+1}(\mathcal{M}(x))$ is also compact. Indeed, $H_{r}^{i+1}$ is outer semicontinuous according to \autoref{lem:PI-ith-feasible},
hence  $H_{r}^{i+1}(\mathcal{M}(x))$ is closed according to \cite[Theorem 5.25(a)]{Rockafellar-Wets-book}. On the other hand, $H_{r}^{i+1}\big(\mathcal{M}(x)\big)$ is bounded since $\mathcal{M}(x)$ itself is bounded and $H_{r}^{i+1}$ is locally bounded, according to \cite[Proposition 5.15]{Rockafellar-Wets-book}. Hence, {$\mathcal{M}(x)\times H_{r}^{i+1}\big(\mathcal{M}(x)\big)$ is compact} and, therefore, the continuity of  $\sigma$ and $\ell$  ensured by \autoref{SA:detect-control} and \autoref{ass:feasibility} is uniform {by Heine theorem}.} {As a result}, for $\frac{\varepsilon}{2\overline N}>0$, there exists $\delta_1>0$ such that for any $\chi_1,\chi_2\in\mathcal{M}(x)\times H_{r}^{i+1}\big(\mathcal{M}(x)\big)$,
\begin{equation}
\lvert\chi_1-\chi_2\rvert\leq\delta_1\Rightarrow\lvert\ell(\chi_1)-\ell(\chi_2)\rvert\leq\tfrac{\varepsilon}{2\overline N}.
\label{eq:elleps}
\end{equation}
On the other hand, for $\varepsilon_1>0$, which is defined in \eqref{eq:eps2}, there exists $\delta_2>0$ such that for any $(x_1,x_2)\in\mathcal{M}(x)\times(\mathcal{M}(x)+\delta_2\mathbb{B})$,
\begin{equation}
\displaystyle\lvert x_1-x_2\rvert\leq\delta_2\Rightarrow\lvert\sigma(x_1)-\sigma(x_2)\rvert\leq\varepsilon_1.\label{eq:sigmaeps2}
\end{equation}
Let $\chi:=(x,u)$ with $u\in H_{r}^{i+1}(x)$. {In view of \autoref{prop:osc-Gi+1}, $G^{i+1}$ is outer semicontinous and locally bounded, it follows that system \eqref{eq:Gi+1} is (nominally) well-posed according to \cite[Theorem 6.30, Assumption 6.5 (A3)]{Goebel-Sanfelice-Teel-book}. Consequently}, the  conditions of \cite[Proposition 6.14]{Goebel-Sanfelice-Teel-book} are verified, {we thus apply this result with the triple $K=\lbrace \chi \rbrace$, $\tau=\overline N$, $\varepsilon=\min\lbrace \delta_1,\delta_2\rbrace$. As a result, there exists $\delta_3>0$ such that for any optimal solution $\psi_{i+1,n}=(\phi_{i+1,n},\upsilon_{i+1,n})\in\mathcal{S}^{i+1}(x_n)$ with $n$ sufficiently large such that $\lvert x-x_n\rvert\leq\delta_3$, there exists a  solution $\hat{\psi}_{i+1,n}=(\hat{\phi}_{i+1,n},\hat{\upsilon}_{i+1,n})\in\mathcal{S}^{i+1}(x)$, not necessarily optimal, such that $\hat{\psi}_{i+1,n}$ and $\psi_{i+1,n}$ are $(\tau,\varepsilon)$-close, see \cite[Definition 5.23]{Goebel-Sanfelice-Teel-book}, i.e., for all $k\in\lbrace 0,\ldots,\overline N\rbrace$},
\begin{equation}
\lvert \psi_{i+1,n}(k,x_n) - \hat{\psi}_{i+1,n}(k,x)\rvert\leq\min\lbrace\delta_1,\delta_2\rbrace.
\label{eq:phihi-and-phihat}
\end{equation}
Given $\psi_{i+1,n}$ and $\hat{\psi}_{i+1,n}$, we now ensure that they lie in $\mathcal{M}(x)\times H_{r}^{i+1}(\mathcal{M}(x))$ for $k\in\lbrace 0,\ldots,\overline N\rbrace$ in order to exploit \eqref{eq:elleps}. For $k\in\{0,\ldots,\overline N\}$, as $\hat{\psi}_{i+1,n}$ is optimal and in view of \eqref{eq:phihi<beta}, $\phi_{i+1,n}(k,x_n)\in\mathcal{M}(x)$ and it follows that $H_{r}^{i+1}(\phi_{i+1,n}(k,x_n))\subseteq H_{r}^{i+1}(\mathcal{M}(x))$. Therefore, for $k\in\lbrace 0,\ldots,\overline N\rbrace$, in view of  \eqref{eq:phihi-and-phihat},  $\hat{\phi}_{i+1,n}(k,x)\in\mathcal{M}(x)+\delta_2\mathbb{B}$, and from \eqref{eq:sigmaeps2} we deduce that
\begin{equation}
\lvert\sigma(\phi_{i+1,n}(k,x_n))-\sigma(\hat{\phi}_{i+1,n}(k,x))\rvert\leq\varepsilon_1.
\label{eq:sigma et psi}
\end{equation}
In view of \eqref{eq:eps2}, $\varepsilon_1\in(0,1]$ and from \eqref{eq:phihi<beta} and  \eqref{eq:sigma et psi}, it follows that \iftoggle{onecolumnversion}{\begin{equation}
\begin{aligned}
\sigma(\hat{\phi}_{i+1,n}(k,x))&\leq\sigma\big(\phi_{i+1,n}(k,x_n)\big)+\varepsilon_1\\
&\leq\beta\big(\sigma(x)+1,0\big)+\varepsilon_1\\
&\leq\beta\big(\sigma(x)+1,0\big)+1.
\end{aligned}\label{eq:ineq:phi_n}
\end{equation}}{$\sigma(\hat{\phi}_{i+1,n}(k,x))\leq\sigma\big(\phi_{i+1,n}(k,x_n)\big)+\varepsilon_1\leq\beta\big(\sigma(x)+1,0\big)+\varepsilon_1\leq\beta\big(\sigma(x)+1,0\big)+1$. }%
We derive that, for all $k\in\lbrace 0,\ldots,\overline N\rbrace$, $\hat{\phi}_{i+1,n}(k,x)$ belongs to $\mathcal{M}(x) $ in view of \eqref{eq:M} and $H_{r}^{i+1}\big(\hat{\phi}_{i+1,n}(k,x)\big)\subseteq H_{r}^{i+1}\big(\mathcal{M}(x)\big)$. Hence, from \eqref{eq:elleps}, we have for all $k\in\lbrace 0,\ldots,\overline N\rbrace$,
\begin{equation}
\lvert\ell\big(\psi_{i+1,n}(k,x_n)\big)-\ell\big(\hat{\psi}_{i+1,n}(k,x)\big)\rvert\leq\frac{\varepsilon}{2\overline N}.\label{eq:ell-psihi-psin}
\end{equation}
We now define $\overline{\psi}_{i+1}(k,x)=\hat{\psi}_{i+1,n}(k,x)$ for $k\in\{0,\ldots,\overline{N}-1\}$ and optimal
     $\overline{\psi}_{i+1}(k,x)=\psi_{i+1}(k-\overline{N},\hat{\phi}_{i+1,n}(\overline{N},x))$ for  $k\in\{\overline{N},\ldots\}$.
We note that $\overline{\psi}_{i+1}\in\mathcal{S}^{i+1}(x)$. As a result, from \eqref{eq:Vi+1},
\begin{equation}\textstyle
\begin{split}
    \MoveEqLeft V_r^{i+1}(x)\\ &\leq\sum\nolimits_{k=0}^\infty \ell(\overline{\psi}_{i+1}(k,x)) \\ & =   \sum\nolimits_{k=0}^{\overline N-1}  \ell(\hat{\psi}_{i+1,n}(k,x)) + V_r^{i+1}(\hat{\phi}_{i+1,n}(\overline N,x)).
    \label{eq:psibar-cost}
\end{split}
\end{equation} 
From \eqref{eq:ell-psihi-psin}, we deduce that
\begin{equation}
V_r^{i+1}(x_n)\geq \sum_{k=0}^{\overline{N}-1} \ell(\psi_{i+1,n}(k,x_n))\geq\sum_{k=0}^{\overline{N}-1} \ell(\hat{\psi}_{i+1,n}(k,x))-\frac{\varepsilon}{2}.\label{eq:sumell-psijpsijhat}
\end{equation}
By adding and subtracting $V_r^{i+1}(\hat{\phi}_{i+1,n}(\overline N,x))$ to \eqref{eq:sumell-psijpsijhat}, in view of \eqref{eq:psibar-cost},{
\begin{align}
V_r^{i+1}(x_n)&\geq \sum_{k=0}^\infty \ell\big(\overline{\psi}_{i+1}(k,x)\big) -\frac{\varepsilon}{2}-V_r^{i+1}(\hat{\phi}_{i+1,n}(\overline N,x))\nonumber\\
&\geq V_r^{i+1}(x) -\frac{\varepsilon}{2}-V_r^{i+1}(\hat{\phi}_{i+1,n}(\overline N,x)).
\label{eq:Vi(xn)>-eps/2+somm 2}
\end{align}
The next inequality comes from $V_r^{i+1}(z)\leq\alphav(\sigma(z))$ for any $z\in\R^{n_x}$ by \autoref{prop:Vi+1Vi} as well as \eqref{eq:sigma et psi},
\begin{align}
\MoveEqLeft V_r^{i+1}\big(\hat{\phi}_{i+1,n}(\overline N,x)\big) \nonumber \\ &{}\leq\alphav\big(\sigma(\hat{\phi}_{i+1,n}(\overline N,x))\big)\leq\alphav\big(\sigma(\phi_{i+1,n}(\overline N,x_n))+\varepsilon_1\big).  \label{eq:Vileqav}
\end{align}
Since $\alphav\in\mathcal{K}_{\infty}$, then for all $a,b\geq0$, $\alphav(a+b)\leq\alphav(2a)+\alphav(2b)$.
{Given $n\in\Zo$} sufficiently big so that \eqref{eq:phihi<beta} holds, it follows from \eqref{eq:eps2}, \eqref{eq:Vileqav} that
$V_r^{i+1}\big(\hat{\phi}_{i+1,n}(\overline N,x)\big) \leq  \alphav\big(2\sigma(\phi_{i+1,n}(\overline N,x_n))\big)+\alphav\big(2\varepsilon_1\big)  \leq   \alphav\big(2\beta(\sigma(x_n),\overline N)\big)+\alphav\big(2\varepsilon_1\big) \leq   \alphav\big(2\beta(\sigma(x)+1,\overline N)\big)+\alphav\big(2\varepsilon_1\big)
 \leq   2\alphav\big(2\varepsilon_1\big)$.
By definition of $\varepsilon_1$ in \eqref{eq:eps2} and the result just above, we deduce
$V_r^{i+1}\big(\hat{\phi}_{i+1,n}(\overline N,x)\big) \leq  \frac{\varepsilon}{2}$. Finally, the combination of the previous inequality and \eqref{eq:Vi(xn)>-eps/2+somm 2}   leads to $V_r^{i+1}(x_n)\geq V_r^{i+1}(x)-\varepsilon$.
We have proved that $V_r^{i+1}$ is lower semicontinuous at $x$, it is thus lower semicontinuous on $\R^{n_x}$ as $x$ has been arbitrarily selected.}

\subsection{Sketch of proof of \autoref{cor:expresult}} \label{appendix:exponential}
Given \autoref{thm:rstability}, it suffices to show that  \autoref{thm:rstability} holds with $\beta:(s,k)\mapsto\dst\tfrac{\overline a_Y}{\underline a_Y}(1-\widetilde{a}_Y)^k s$, which is ${\exp}{-}\KL$. This is the case in view of the conditions of \autoref{cor:expresult}. To see this, note the following. First, 
\autoref{lem:lyapunov-model-based-i+1} is verified with $\rho_V=\mathbb{I
}$ and $\rho_W=\frac{1}{c_W}\mathbb{I
}$ with the functions given in the case where $\chi_W\leq\mathbb{I}$ in the proof of \autoref{lem:lyapunov-model-based-i+1}. Then, by following the steps in the proof of \autoref{lem:KL-stability-i=0}, we derive $\overline{\beta}(s,k)=(1-\widetilde{a}_Y)^{k}s$ such that $Y^{i}(\phi^i(k,x))\leq\overline{\beta}(Y^{i}(x),k)$ for all $x\in\R^{n_x}$. The desired result is derived by using the conditions of \autoref{cor:expresult}.

%%%%%%%%%%%%%%%%%%%%%%%%%%%%%

\vspace{-1ex}

\section{\hspace{2em} Proofs of \autoref{subsec:near-opti}}\label{appendix:secnear-opti}
\subsection{Proof of \autoref{thm:near-opti}}\label{appendix:near-opti}

Let $x\in\R^{n_x}$ and $i\in\Zp$. From item (i) of \autoref{prop:Vi+1Vi} and   Bellman equation,  it follows $V_r^i(x)-V^\star(x)
=\ell(x,h_r^{\star,i}(x))+ V_r^{i}(\phi(1,x,h_r^{\star,i}))
    -\ell(x,h^\star(x))-V^\star(\phi(1,x,h^\star)),$
for any $h_r^{\star,i}\in H_r^{\star,i}$ and $h^{\star}\in H^{\star}$. Moreover,  $V_r^i(x)=\ell(x,h_r^{\star,i}(x))+V_r^{i}(\phi(1,x,h_r^{\star,i}))\leq \ell(x,h^{i}(x))+V_r^{i}(f(x,h^i(x)))\leq\ell(x,u)+V_r^{i-1}(f(x,u))$ for any $h^i\in H^i$ and $u\in\mathcal{U}(x)$ as consequence of \eqref{eq:PI-ith-improvement}, \eqref{eq:PI-ith-improvement-best} and item (ii) of \autoref{prop:Vi+1Vi}. We derive $V_r^i(x)-V^\star(x)\leq\ell(x,u)+ V_r^{i-1}(f(x,u))-\ell(x,h^\star(x))-  V^\star(\phi(1,x,h^\star))\leq\ell(x,h^\star(x))+V_r^{i-1}(\phi(1,x,h^\star))-\ell(x,h^\star(x))- V^\star(\phi(1,x,h^\star))= V_r^{i-1}(\phi(1,x,h^\star)) -V^\star(\phi(1,x,h^\star))$.
We repeat the above reasoning $i-1$ times, hence $V_r^i(x)-V^\star(x) \leq V_r^{0}(\phi(i,x,h^\star)) -V^\star(\phi(i,x,h^\star))$.
The first part of \autoref{thm:near-opti} is obtained. 

We now show \eqref{eq:explicit-near} holds  and we distinguish two cases as in the proof of \autoref{lem:lyapunov-model-based-i=0} in \autoref{appendix:profLyapunov}.

\textbf{Case where $\bm{\chi_W\leq\mathbb{I}}$:} For any $z\in\R^{n_x}$, $V_r^0(z)-V^\star(z) = V_r^0(z)-V^\star(z)\pm W(z)=Y^0(z)-Y^\star(z)$, where $Y^\star=V^\star+W$ defined in \autoref{sssec:robustpractical}, in view of the proof of  \autoref{lem:KL-stability-i=0} in \autoref{appendix:profLyapunov}. Moreover that $V_r^0(z)-V^\star(z)\leq V_r^0 (z)$ since $V^\star(z)\geq0$. Hence $V_r^0(z)-V^\star(z)\leq\min\{V_r^0(z),Y^0(z)-Y^\star(z)\}$. On the one hand, $V_r^{0}(z)\leq\overline{\alpha}_V(\sigma(z))$  for any $z\in\R^{n_x}$ from \autoref{SA:stabilizing-initial-policy}, and on the other hand,
$\underline{\alpha}_Y(\sigma(z))\leq Y^\star(z)$ moreover $Y^0(z)\leq \overline{\alpha}_Y(\sigma(z))$, $Y^0(z)-Y^\star(z)\leq\overline{\alpha}_Y(\sigma(z))-\underline{\alpha}_Y(\sigma(z))$.  Therefore  $V_r^0(z)-V^\star(z)\leq\min\{V_r^0(z),Y^0(z)-Y^\star(z)\} \leq \min\{\overline{\alpha}_V(\sigma(z)),\overline{\alpha}_Y(\sigma(z))-\underline{\alpha}_Y(\sigma(z))\}=: \widehat{\alpha}(\sigma(z))$, where $\widehat{\alpha} \circ \sigma$ is continuous and positive semidefinite but not necessarily non-decreasing.
Thus, we construct $\widetilde{\alpha}(s):= \max_{\hat{s}\in[0,s]} \widehat{\alpha}(\hat{s})$  and we obtain suitable $\widetilde{\alpha}$ continuous, positive definite and non-decreasing such that $V_r^0(z)-V^\star(z)\leq\widetilde{\alpha}(\sigma(z))$, and this concludes the case where $\chi_W\leq\mathbb{I}$.

\textbf{Case where $\bm{\chi(s)>s}$ for some $\bm{s>0}$:} For any $z\in\R^{n_x}$,
since $V^\star(z)\geq \ell(z,u)$  for any $u\in\mathcal{U}(x)$ and $V_r^{0}(z)\leq\overline{\alpha}_V(\sigma(z)))$  for any $z\in\R^{n_x}$ from \autoref{SA:stabilizing-initial-policy},  $V_r^0(z)-V^\star(z) \leq \overline{\alpha}_V(\sigma(z))-\ell(z,u)$. From \autoref{SA:detect-control} and $\ell(z,u)\geq0$ for any $(z,u)\in\mathcal{W}$, we have
$\max\{0,\alpha_W(\sigma(z))-\overline{\alpha}_W(\sigma(z))\} \leq\chi_W(\ell(z,u))$ hence
$-\ell(z,u)\leq-\chi_W^{-1}\left(\max\{0,\alpha_W(\sigma(z))-\overline{\alpha}_W(\sigma(z))\}\right)=:-\widetilde{\alpha}_W(\sigma(z))$, where $\widetilde{\alpha}_W \circ \sigma$ is continuous and positive semidefinite but not necessarily class $\K_\infty$. Nevertheless, $V_r^0(z)-V^\star(z) \leq \overline{\alpha}_V(\sigma(z))-\ell(z,u)\leq \overline{\alpha}_V(\sigma(z))-\widetilde{\alpha}_W(\sigma(z))=:\widehat{\alpha}(\sigma(z))$.
As done above,  we construct $\widetilde{\alpha}(s):= \max_{\hat{s}\in[0,s]} \widehat{\alpha}(\hat{s})$  and we obtain suitable $\widetilde{\alpha}$ continuous, positive definite and non-decreasing such that $V_r^0(z)-V^\star(z)\leq\widetilde{\alpha}(\sigma(z))$, and  this concludes the  case where $\chi(s)>s$ for some $s>0$.

All that remains is to show that $\sigma(\phi(i,x,h^\star))\leq\beta(\sigma(x),i)$ holds for any $x\in\R^{n_x}$ and $i\in\Zo$. To see this, it suffices to note that, in view of \autoref{SA:well-posed} and \autoref{SA:stabilizing-initial-policy}, $V^\star\leq V^0\leq \overline{\alpha}\circ\sigma$, and then to follow the same steps as in the proof of \autoref{lem:lyapunov-model-based-i=0} in \autoref{appendix:profLyapunov} and the proof of \autoref{lem:KL-stability-i=0} in \autoref{appendix:proofKL-stability} for $h^\star\in H^\star$ as in \eqref{eq:Hstar} instead of $H_r^{\star,i+1}$.  Then,  since $\widetilde{\alpha}$ is non-decreasing, $V_r^i(x)-V^\star(x) \leq \widetilde{\alpha}(\beta(\sigma(x),i))$
and this concludes the proof. 

%%%%%%%%%%%%

\subsection{Proof of \autoref{prop:uniform}} \label{appendix:proofuniform}

Item (i) of \autoref{prop:uniform} holds in view of item (ii) of \autoref{prop:Vi+1Vi}. We now focus on item (ii) of \autoref{prop:uniform}.
Let $\delta>0$ and $K\subset\R^{n_x}$ compact. Let $\Delta>0$ sufficiently large  such that $\{z\in\R^{n_x}\,:\, \sigma(z)\leq\Delta\}\supseteq K$ and  let $x\in K$, thus $\sigma(x)\leq\Delta$. Since $\beta\in\KL$ and $\widetilde{\alpha}$ is continuous and zero at zero, there exists $i^\star$ sufficiently large that
$\widetilde{\alpha}(\beta(\Delta,i^\star))\leq\delta$.
As $\widetilde{\alpha}$ is non-decreasing and $\beta\in\KL$,  $\widetilde{\alpha}(\beta(\sigma(x),i))\leq\widetilde{\alpha}(\beta(\Delta,i))\leq\widetilde{\alpha}(\beta(\Delta,i^\star))\leq\delta$ for $i\geq i^\star$, and the proof is finished by invoking \eqref{eq:explicit-near}.
%%%%%%%%%%%%%%%%%%%%%%%%%%%%%

\section{\hspace{2em} Sketch of proofs of \autoref{sec:sys-cost-problem}}  \label{appendix:PI}

In view of \autoref{ass:rfeasibility}, we do not require the interplay between stability and feasibility as in \autoref{prop:induction-statement} in \refappendix{appendix:recursive}. Hence, to derive stability properties as in \autoref{thm:rstability-PI}, it suffices to show that: (i) for all $i\in\Zo$, $x\in\R^{n_x}$ and  any $h^{i+1}\in H^{i+1}$, $J(x,h^{i+1})\leq V^i(x)$, which is standard in PI literature, see e.g., \cite[Section 4.2]{Sutton-Barto-book2018}, and is akin to \autoref{prop:Vi+1Vi} for \PIplus; (ii)  invoke
$V^0(x)\leq\overline{\alpha}_V \circ \sigma(x)$ for any $x\in\R^{n_x}$ by \autoref{SA:stabilizing-initial-policy} and then to follow the Lyapunov arguments made in the proof of \hyperref[{lem:lyapunov-model-based-i=0}]{Lemmas \ref*{lem:lyapunov-model-based-i=0}} \hyperref[lem:KL-stability-i=0]{and \ref*{lem:KL-stability-i=0}}. On the other hand, for near-optimality of PI in \autoref{thm:near-opti-PI} and \autoref{prop:uniform-PI}, it suffices to follow the proofs in \refappendix{appendix:secnear-opti} with $V^i$ in place of $V_r^i$.

\vspace{-1ex}

\section{\hspace{2em} Proofs of \autoref{sec:robust}} \label{appendix:robustness}

\subsection{Proof of \autoref{thm:robust-practical-stability}} \label{appendix:proofrobust-practical-stability}

As explained in \autoref{sssec:robustdefinition}, we will invoke \cite[Theorem 2.8]{Kellett-Teel-siam-jco-05} to obtain robust $\KL{-}\text{stability}$ with respect to some measures $(\sigma_1,\sigma_2)$ on $\mathcal{X}\subseteq \R^{n_x}$. The first requirement, regarding the compactness of $F^i_r(x)$ for all $x\in\mathcal{X}$ is established in \autoref{lem:compactness}. All that remains is  to show that $Y^\star=\rho_V(V^\star)+\rho_W(W)$ is continuous and is a suitable Lyapunov function for the considered  stability  property of system \eqref{eq:ith-auto-sys} with respect to measures $(\sigma_1,\sigma_2)$ on $\mathcal{X}\subseteq \R^{n_x}$.

First, we show in the next lemma that $V^\star$ satisfies a useful dissipation inequality along \eqref{eq:ith-auto-sys}.
\begin{lem} \label{lem:practicaldecrease}For any $i\in\Zo$, $x\in\R^{n_x}$, $\upsilon\in F_r^{i}(x)$, $h_r^{\star,i}\in H_r^{\star,i}$ such that $f(x,h^{\star,i}(x))=\upsilon$ and $h^{\star}\in H^{\star}$, $V^\star(\upsilon)-V^\star(x)\leq -\ell(x,h_r^{\star,i}(x)) + (V_r^0-V^\star)(\phi(i,x,h^\star))$.
 \mbox{} \hfill $\Box$ 
\end{lem}
\textbf{Proof:} Let $i\in\Zo$, $x\in\R^{n_x}$, $\upsilon\in F_r^{i}(x)$  and $h^{\star,i}\in H^{\star,i}$ such that $f(x,h^{\star,i}(x))=\upsilon$. Since $V_r^i(x)=\ell(x,h^{\star,i}(x))+V_r^i(\upsilon)$ and $V^\star\leq V_r^i$ for all $i\in\Zp$ by definition of $V^\star$ in \eqref{eq:Vstar}, $V_r^i(x)\geq\ell(x,h^{\star,i}(x))+V^\star(\upsilon)$ holds. In view of \eqref{eq:abstract-near},  $V_r^i(x)-V^\star(x)\leq (V_r^0-V^\star)(\phi(i,x,h^\star))$, and the desired result holds by combining the previous inequalities.
\mbox{}\hfill $\blacksquare$

The details of the remaining of the proof is omitted for space reasons. The proof  follows by employing similar steps  as done in the proof of \autoref{lem:lyapunov-model-based-i+1} to derive there exists $i^\star\in\Zo$ such that, for any $i>i^\star$,
$\sigma_1(\phi^i(k,x))= \sigma(\phi^i(k,x))-\delta\leq \max\{\widehat{\beta}(k,\sigma(x))-\delta,0\}\leq\widehat{\beta}(k,\sigma_2(x))$ holds for any $x\in\mathcal{X}:=\{z\in\R^{n_x}\,:\, \sigma(z)<\Delta\}$. The Lyapunov function used to establish this stability property is $Y^\star$, which is continuous  in view  of \autoref{ass:Ucont} and following similar steps as \cite[Theorem 6]{Kellett-Teel-siam-jco-05}. As $F_r^i$ is non-empty and compact on $\mathcal{X}$ in view of \autoref{lem:compactness}, we invoke \cite[Theorem 2.8]{Kellett-Teel-siam-jco-05} to establish the desired result.
 
%%%%%%%%%%%%%%

\subsection{Proof of \autoref{prop:continuous}} \label{appendix:continuous}

We first show that, under \autoref{ass:samecost}, when  $V_r^i$ is continuous for some $i\in\Zo$, $V_r^{i+1}$ is continuous. Since $V_r^0$ is continuous, the desired result  then follows by induction.

\begin{lem} \label{prop:cont}
When \autoref{ass:samecost} holds, $V_r^i$  continuous on $\R^{n_x}$  for $i\in\Zo$ implies $V_r^{i+1}$ continuous on $\R^{n_x}$. \mbox{}\hfill $\Box$
\end{lem}
\noindent\textbf{Proof:} In light of \autoref{prop:semilowercontinuity}, all that remains to be shown is that  $V_r^{i+1}$ is upper semicontinuous on $\R^{n_x}$ given $V_r^i$ continuous on $\R^{n_x}$ and \autoref{ass:samecost}. Hence, we extend here the proof of \autoref{prop:semilowercontinuity} in \autoref{appendix:lowersemicontinuity}  to show, for any $\varepsilon>0$, for any $x\in\R^{n_x}$ and sequence $x_n\to x$, there exists $N\in\Zp$, for any $n\geq N$ such that $V^{i+1}(x)\geq V^{i+1}(x_n)-\varepsilon$. 

Let $x\in\R^{n_x}$, $\varepsilon\in\Rlo$ and $n$ sufficiently large such that any invoked inequalities of the proof of \autoref{prop:semilowercontinuity} in \autoref{appendix:lowersemicontinuity} hold in the sequel.  Given  the continuity of $V_r^i$, we have that $H^{i+1}$ is outer semicontinuous in view of \cite[Theorem 7.41(b)]{Rockafellar-Wets-book}, hence $H^{i+1}=H_r^{i+1}$, and in view of \autoref{ass:samecost} we have $H^{i+1}=H_r^{i+1}=H_r^{\star,i+1}$. Thus, for any two solutions $\psi,\widetilde{\psi}\in\mathcal{S}^{i+1}(z)$, $\sum_{k=0}^{\infty} \ell(\psi(k,z))=\sum_{k=0}^{\infty} \ell(\widetilde{\psi}(k,z))=V_r^{i+1}(z)$ for all $z\in\R^{n_x}$, that is, all solutions are optimal. Therefore,  we have $V^{i+1}(x)=\sum_{k=0}^{\infty} \ell(\hat{\psi}_{i+1,n}(k,x))$ for $\hat{\psi}_{i+1,n}=(\hat{\phi}_{i+1,n},\hat{v}_{i+1,n})\in\mathcal{S}^{i+1}(x)$ as defined in \autoref{appendix:lowersemicontinuity}. From \eqref{eq:elleps}, we deduce that $V^{i+1}(x)\geq \sum_{k=0}^{\overline{N}-1} \ell(\hat{\psi}_{i+1,n}(k,x))\geq\sum_{k=0}^{\overline{N}-1} \ell({\psi}_{i+1,n}(k,x_n))-\frac{\varepsilon}{2}$ for $\psi_{i+1,n}=(\phi_{i+1,n},v_{i+1,n})\in\mathcal{S}^{i+1}(x_n)$  as defined in \autoref{appendix:lowersemicontinuity}.
By adding and subtracting $\sum_{k=\overline N}^\infty \ell({\psi}_{i+1,n}(k,x_n))$ to the previous inequality and in view of the optimality of ${\psi}_{i+1,n}$, $V^{i+1}(x)\geq \sum_{k=0}^\infty \ell({\psi}_{i+1,n}(k,x_n))-\frac{\varepsilon}{2}-\sum_{k=\overline{N}}^\infty \ell({\psi}_{i+1,n}(k,x_n))= V^{i+1}(x_n)-\frac{\varepsilon}{2}-\sum_{k=\overline{N}}^\infty \ell({\psi}_{i+1,n}(k,x_n))$.
Again, in view of the optimality of $\psi_{i+1,n}$, $V^{i+1}({\phi}_{i+1,n}(\overline{N},x_n))=\sum_{k=\overline N}^\infty \ell\big({\psi}_{i+1,n}(k,x_n)\big)$. Hence it follows from the previous inequality that
\begin{equation}
V^{i+1}(x)\geq V^{i+1}(x_n)-\frac{\varepsilon}{2}-V^{i+1}({\phi}_{i+1,n}(\overline{N},x_n)).\label{eq:bis34}
\end{equation}
The following inequality is deduced from $V^{i+1}(z)\leq\alphav(\sigma(z))$ for any $z\in\R^{n_x}$ in view of \autoref{prop:Vi+1Vi}, \autoref{SA:stabilizing-initial-policy}, \eqref{eq:sigmaeps2}, and optimality of $\psi_{i+1,n}$,
\iftoggle{onecolumnversion}{\begin{equation}
\begin{split}
    V^{i+1}\big({\phi}_{i+1,n}(\overline N,x_n)\big)&{}\leq\alphav\big(\sigma({\phi}_{i+1,n}(\overline N,x_n))\big)\\&{}\leq\alphav\big(\sigma(\hat{\phi}_{i+1,n}(\overline N,x))+\varepsilon_1\big).
\end{split}
\label{eq-Vileqavupper}
\end{equation}}{$V^{i+1}\big({\phi}_{i+1,n}(\overline N,x_n)\big)\leq\alphav\big(\sigma({\phi}_{i+1,n}(\overline N,x_n))\big)\leq\alphav\big(\sigma(\hat{\phi}_{i+1,n}(\overline N,x))+\varepsilon_1\big)$. }%
{Given $n\in\Zo$} sufficiently big so that \eqref{eq:phihi<beta} holds, by proceeding like in the end of the proof of \autoref{prop:semilowercontinuity}, we obtain 
$V^{i+1}\big({\phi}_{i+1,n}(\overline N,x_n)\big)\leq2\alphav\big(2\varepsilon_1\big)$.
By definition of $\varepsilon_1$ in \eqref{eq:eps2} we deduce $V^{i+1}\big({\phi}_{i+1,n}(\overline N,x_n)\big) \leq  \frac{\varepsilon}{2}$.
Finally, in view of \eqref{eq:bis34}, $V^{i+1}(x)\geq V^{i+1}(x_n)-\varepsilon$.
Hence, $V^{i+1}$ is upper semicontinuous at $x$.  We have proved that $V^{i+1}$ is continuous at $x$, it is thus continuous on $\R^{n_x}$ as $x$ has been arbitrarily selected.

%%%%%%%%%%%%%%%%%

\subsection{Proof of \autoref{thm:robust-asymptotic-stability}} \label{appendix:proof-robust-asymptotic-stability}

The proof is an application of \cite[Theorem 2.8]{Kellett-Teel-siam-jco-05} as by \autoref{lem:compactness}, \autoref{thm:rstability} and \autoref{prop:continuous} all the required conditions hold.

%%%%%%%%%%%%
\vspace{-2ex}

\section{\hspace{2em} Technical lemmas} \label{appendix:technical}

Three technical lemmas used in \refappendix{appendix:recursive} are given here.

\begin{lem}\label{appendix:deltaosclsc}  Given an outer-semicontinuous set-valued map $H:\R^{n}\rightrightarrows \R^{m}$ with $n,m\in\Zp$, the function $(x,u)\mapsto\delta_{H(x)}(u)$ is lower semicontinuous on $\R^n\times\R^m$. 
\mbox{}\hfill $\Box$
\end{lem}
\textbf{Proof:} Let $x\in\R^n$ and $u\in\R^m$. To prove that $\delta_{H(\cdot)}(\cdot)$ is lower semicontinuous at $(x,u)$, we will show that $\liminf_{(\tilde x,\tilde u)\to (x,u)}\allowbreak \delta_{H(\tilde x)}(\tilde u)\geq \delta_{H(x)}(u)$, see \cite[Definition 1.5]{Rockafellar-Wets-book}. This is equivalent to $\min\{\alpha\in\overline{\R}\,:\, \exists (x_n,u_n)\to(x,u), \delta_{H(x_n)}(u_n)\to\alpha\}\geq\delta_{H(x)}(u)$, as  $\liminf_{(\tilde x,\tilde u)\to (x,u)} \delta_{H(\tilde x)}(\tilde u) = \min\{\alpha\in \overline{\R}\,:\, \exists (x_n,u_n)\to(x,u), \delta_{H(x_n)}(u_n)\to\alpha\}$ by \cite[Lemma 1.7]{Rockafellar-Wets-book}.  For what follows, we   employ the fact that both $\min\{\alpha\in\overline{\R}\,:\, \exists (x_n,u_n)\to(x,u), \delta_{H(x_n)}(u_n)\to\alpha\}$ and $\delta_{H(x)}(u)$ only take values  in $\{0,\infty\}$.

We distinguish two cases.  First, when $\min\{\alpha\in\overline{\R}\,:\, \exists(x_n,u_n)\to(x,u), \delta_{H(x_n)}(u_n)\to\alpha\}=\infty$. The desired result holds. 
Second, when $\min\{\alpha\in\overline{\R}\,:\, \exists(x_n,u_n)\to(x,u), \delta_{H(x_n)}(u_n)\to\alpha\}=0$. This implies that there exists a (sub)sequence $(x_n,u_n)\to(x,u)$ such that $u_n\in H(x_n)$ as $\delta_{H(x_n)}(u_n)\to0$. Hence, since $H$ is outer semicontinuous, we have $u\in H(x)$. Therefore, $\delta_{H(x)}(u)=0$, and $\min\{\alpha\in\overline{\R}\,:\, \exists(x_n,u_n)\to(x,u), \delta_{H(x_n)}(u_n)\to\alpha\}=0=\delta_{H(x)}(u)$.
We have proved  that $\delta_{H(\cdot)}(\cdot)$ is lower semicontinuous at arbitrary $(x,u)$, it is therefore lower semicontinuous on
$\R^{n}\times\R^m$. \mbox{} \hfill $\blacksquare$

\begin{lem}\label{lem:levelboundeduniform} Let $n,m \in\Zp$ and $g,\ell:\R^n\times\R^m\to\Rlo$ with $g(x,u)\geq \ell(x,u)$ for any $(x,u)\in\R^{n}\times\R^{m}$. If $\ell$ is level-bounded in $u$, locally uniform in $x$,  then so is $g$.
\mbox{}\hfill $\Box$
\end{lem}
\textbf{Proof:} Let $x\in\R^n$, $\alpha\in\R$, $\mathcal{S}$ be a neighbourhood of $x$ and  $B\subset\R^{m}$ be bounded such that $\{u\in\R^m\,:\,\ell(z,u)\leq\alpha\}\subset B$ for any $z\in\mathcal{S}$, which exist by item (ii) of \autoref{ass:feasibility}, see \autoref{def:levelboundeduniform}.  Necessarily $\{u\in\R^m\,:\,g(z,u)\leq\alpha\}\subseteq\{u\in\R^m\,:\,\ell(z,u)\leq\alpha\}$ for  $z\in\mathcal{S}$, as  as $\ell(z,u)\leq g(z,u)$, hence $\{u\in\R^m\,:\,g(z,u)\leq\alpha\}\subset B$  for any $z\in\mathcal{S}$. As $x$ and $\alpha$ have been arbitrarily, the desired result holds.  \mbox{} \hfill $\blacksquare$

\begin{lem}\label{lem:generatingbeta} Let
$\alpha:\Rlo\to\Rlo$ with $\alpha(0)=0$, for any sequence $s_k\in\Rlo$, $k\in\Zo$, with $s_{k+1}\leq \alpha(s_{k})$, we have $s_{k+1}\leq\widetilde{\beta}(s_0,k):=\max_{\hat{s}\in[0,s_0]}\alpha^{(k)}(\hat{s})$. Moreover, when $\alpha$ is continuous and $\alpha(s)< s$ for all $s\in\Rlp$, $\widetilde{\beta}\in\KL$. 
\mbox{}\hfill $\Box$
\end{lem}
\textbf{Proof:} To show $s_k\leq\widetilde{\beta}(k,s_0)$, we  proceed by induction. For $k=0$,  $s_0=\widetilde{\beta}(0,s_0)$.
We assume $s_k\leq\widetilde{\beta}(k,s_0)$ and  show $s_{k+1}\leq\widetilde{\beta}(k+1,s_0)$. On the one hand, $s_{k+1}\leq \alpha(s_k)\leq \max_{\hat s\in[0,s_k]}\alpha(\hat {s})\leq \max_{\hat s\in[0,\widetilde{\beta}(k,s_0)]}\alpha(\hat {s})$  as $s_k\leq\widetilde{\beta}(k,s_0)$ by the induction hypothesis. On the other hand, $\max_{\hat s\in[0,\widetilde{\beta}(k,s_0)]}\alpha(\hat {s})=\max_{\hat{s}\in[\underline{s},\overline{s}]}\alpha(\hat{s})$ with $\underline{s}:=0$ and $\overline{s}:=\widetilde{\beta}(k,s_0)$.  Noting that $ 0=\min_{\hat{\hat{s}}\in[0,s_0]}\alpha^{(k)}(\hat{\hat{s}})$  and $\widetilde{\beta}(k,s_0)=\max_{\hat{\hat{s}}\in[0,s_0]}\alpha^{(k)}(\hat{\hat{s}})$, we have $\max_{\hat{s}\in[\underline{s},\overline{s}]}\alpha(\hat{s})=\max_{\hat{s}\in[0,s_0]}\alpha(\alpha^{(k)}(\hat{s}))$ by the definitions of $\underline{s}$ and $\overline{s}$. We conclude  $s_{k+1}\leq\max_{\hat{s}\in[0,s_0]}\alpha(\alpha^{(k)}(\hat{s}))=\widetilde{\beta}(k+1,s_0)$. 

Now we establish $\widetilde{\beta}\in\KL$ when $\alpha$ is continuous and $\alpha(s)<s$ for all $s>0$. First, it is clear that $\alpha^{(k)}$ is continuous and zero at zero for any $k\in\Zo$, hence $\widetilde{\beta}(\cdot,k)$ is continuous, zero at zero and non-decreasing for any $k\in\Zo$, thus $\widetilde{\beta}$ is class-$\K$ in its first argument. Moreover, as $\alpha(s)<s$ for all $s>0$ and $\alpha(0)=0$, we have
$\widetilde{\beta}(s,k+1)<\widetilde{\beta}(s,k)$ when $\widetilde{\beta}(s,k)>0$ and $\widetilde{\beta}(s,k+1)=\widetilde{\beta}(s,k)=0$ otherwise. Therefore $\widetilde{\beta}(s,\cdot)$ monotonically decreases to zero.  \mbox{} \hfill $\blacksquare$

%\begin{rem} XXX Argument for tightness of beta
%We note that $\widetilde{\beta}$ is equal to $(s,k)\mapsto\overline{\alpha}^{(k)}(s)$ where $\overline{\alpha}=\max_{\hat{s}\in[0,s]}\alpha(\hat{s})$ is non-decreasing. To see this, it suffices to take a subsequence of $s_k$ such that $s_{k+1}=\overline{\alpha}(s_k)$ for $s_1\leq\alpha(s_0)$ and $s_0\in\Rlo$. In this way,  $\overline{\alpha}^{(k)}(s_0)\leq\widetilde{\beta}(s_0,k)$ for any $k\in\Zo$, and moreover that $\widetilde{\beta}(s,k)\leq\overline{\alpha}^{(k)}(s)$ in view of $\widetilde{\beta}(s,k+1)\leq\overline{\alpha}(\widetilde{\beta}(s,k))$ for any $s\in\Rlo$ and $k\in\Zo$. This also shows that  $\widetilde{\beta}$ is the tightest non-decreasing upper-bound, as it is equal to a subsequence of $s_k$. \mathieu{[Issue at $k=1$ when (and only when) $\alpha(s_0)<\overline{\alpha}(s_0)$]}. \mathieu{[This needs an extra push. The math leads me to consider instead $\max_{\hat{s}\in[0,\alpha(s_0)]}\alpha^{(k-1)}(\hat{s})$ and $\overline{\alpha}^{(k-1)}\circ\alpha$ to obtain equality. Note they are \emph{not} non-decreasing: the least non-decreasing upper-bound are still given by $\widetilde{\beta}$ or equivalently $(s,k)\mapsto\overline{\alpha}^{(k)}(s)$!]}\mbox{} \hfill $\Box$
%\end{rem}
\bibliographystyle{plain}
\iftoggle{onecolumnversion}{\bibliography{IEEEabrv,bib_global}}{\bibliography{bib_global_abbrv}}

\vspace{-1cm}

\begin{IEEEbiography}[{\includegraphics[width=1in,height=1.25in,clip,keepaspectratio]{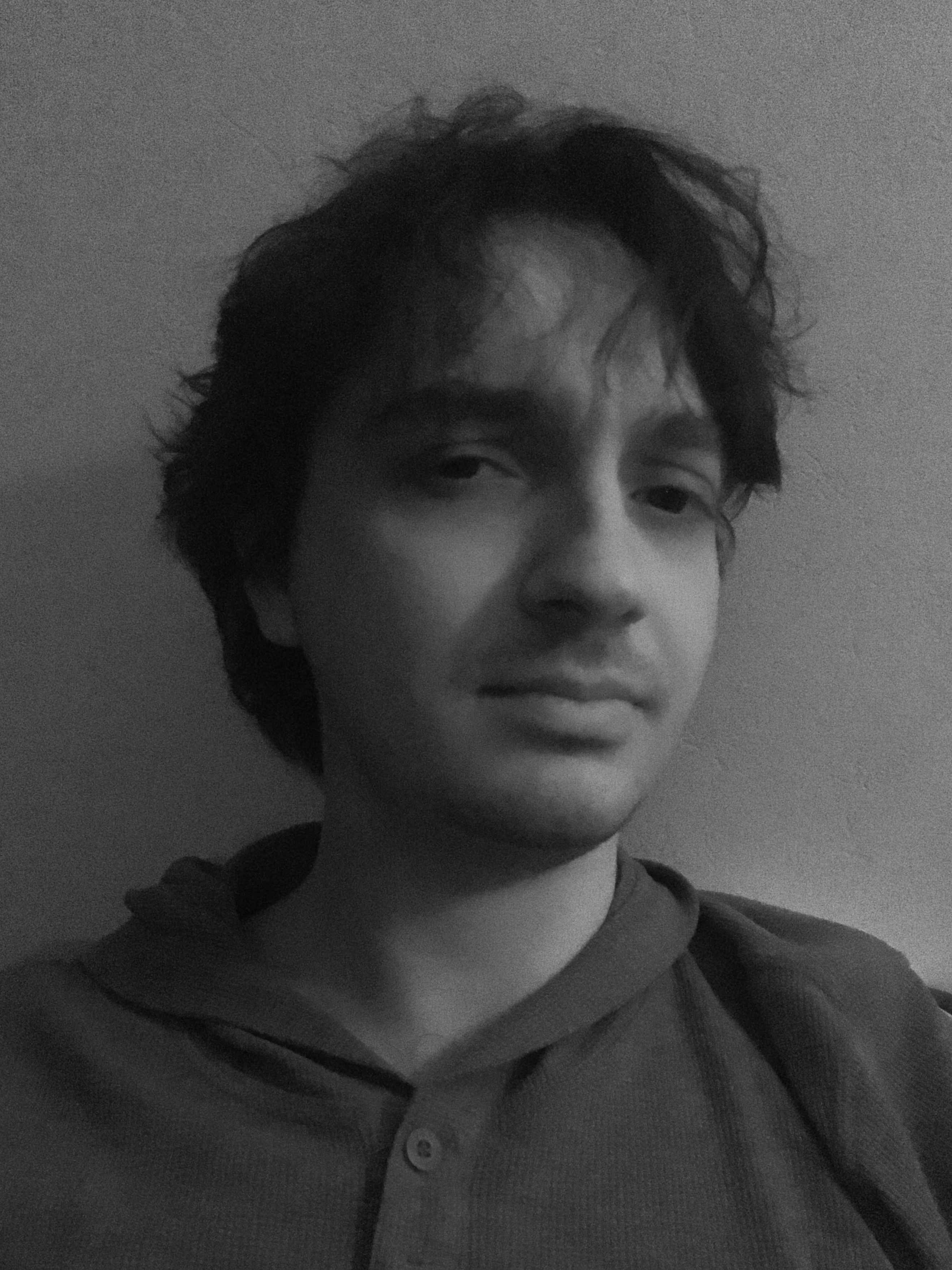}}]{Mathieu Granzotto}
received his engineering degree in Control and Automation in 2016 from ``Universidade Federal de Santa Catarina'', Brazil. In 2019, he received his Ph.D. in Control
Theory from Universit\'e de Lorraine, France, where he was a Temporary Research and Teaching Attaché at CRAN. Since 2022, he is a research fellow at the  Department of Electrical  and Electronic Engineering  (DEEE) at the  University of
Melbourne, Australia. %His research interests include non-linear systems, optimal control and dynamic programming methods.
\vspace{-0.5cm}
\end{IEEEbiography}

\begin{IEEEbiography}[{\includegraphics[width=1in,height=1.25in,clip,trim={0cm 3cm 0 1cm},keepaspectratio]{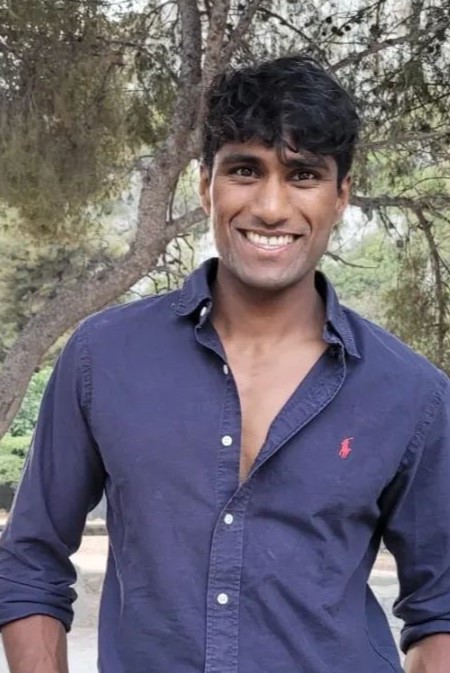}}]{Olivier  Lindamulage De Silva}
received the ``Ingénieur'' degree in Digital Systems Engineering from ENSEM (France) in 2020. He obtained the M.Sc. by Research in Robotic Vision Learning from Université de Lorraine (France) in 2020. He is currently a Ph.D student in Control Theory at Université de Lorraine.
\vspace{-0.5cm}
\end{IEEEbiography}

\begin{IEEEbiography}[{\includegraphics[width=1in,height=1.25in,clip,keepaspectratio]{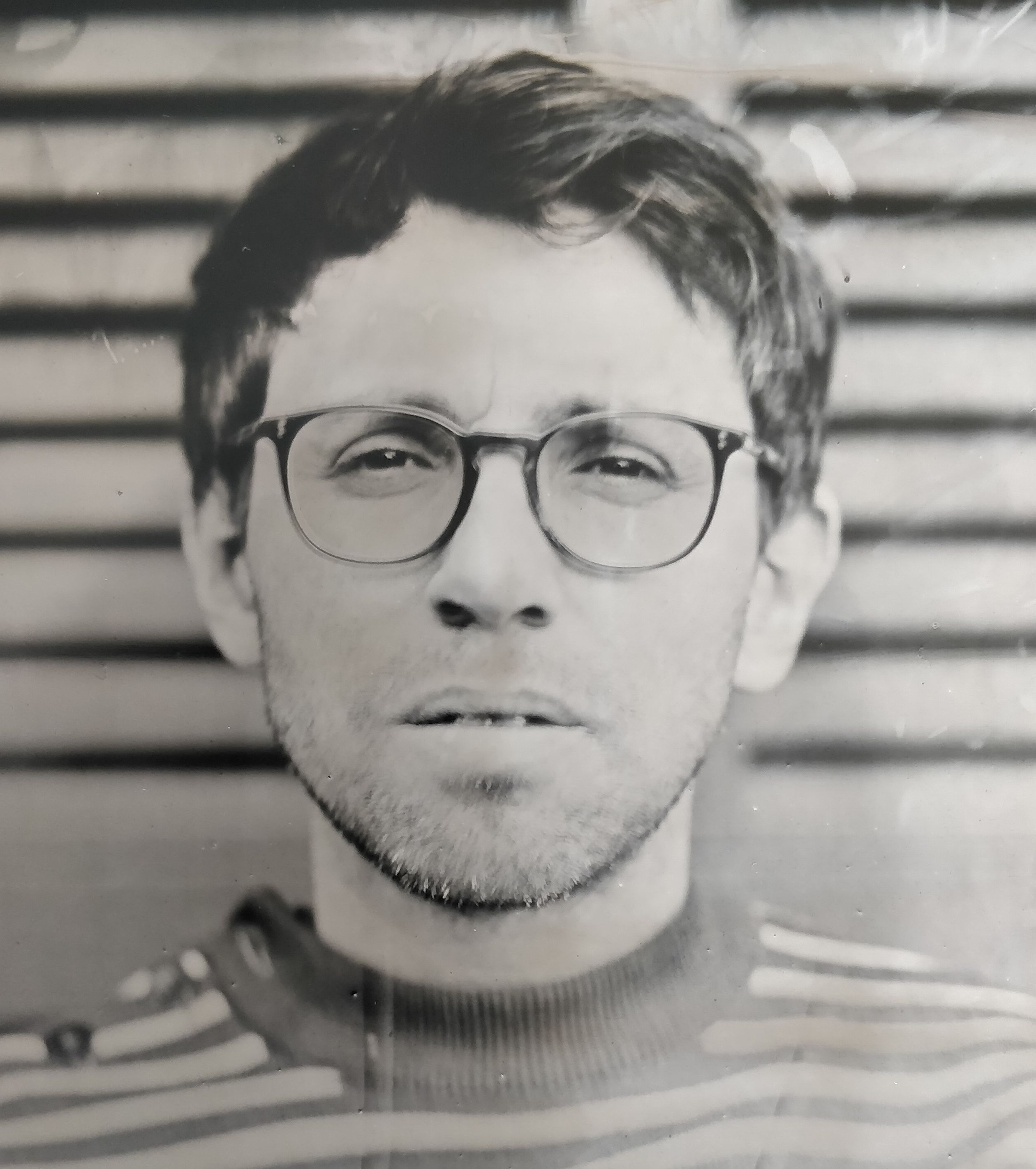}}]{Romain Postoyan}
received the ``Ingénieur''  degree in Electrical and Control Engineering  from ENSEEIHT (France) in 2005.  He obtained the M.Sc.
by Research in Control Theory \& Application from Coventry  University (United Kingdom) in 2006 and the Ph.D. in Control
Theory from Universit\'e Paris-Sud (France) in 2009.  In  2010, he  was a  research assistant at  the University  of Melbourne (Australia).  Since 2011,  he is  a CNRS
researcher at the CRAN (France). %He obtained the ``Habilitation à Diriger des Recherches'' from Université de Lorraine in 2019. %He served/serves as an Associate Editor for the
%journals:  Automatica,  IEEE  Control  Systems  Letters  and  IMA  Journal  of  Mathematical  Control  and  Information.
\vspace{-0.5cm}
\end{IEEEbiography}

\begin{IEEEbiography}[{\includegraphics[width=1in,height=1.25in,clip,trim={0cm 0cm 0 1cm},keepaspectratio]{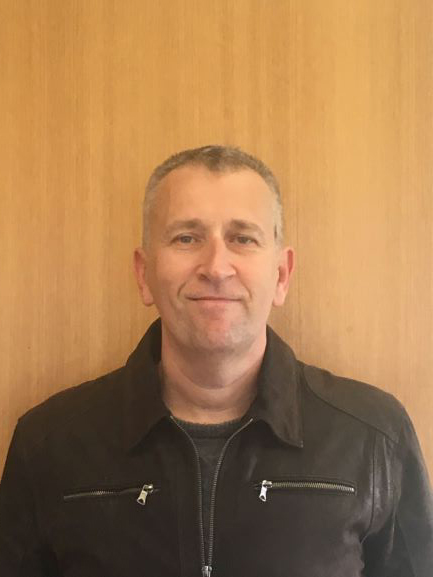}}]{Dragan Nešić}
is a Professor at  The University of Melbourne, Australia. He currently serves as Associate Dean Research at the Melbourne School of Engineering.  His research interests include networked control systems, reset systems, extremum seeking control, hybrid control systems, event-triggered control, security and privacy in cyber-physical systems, and so on. 
He is a Fellow of the Institute of Electrical and Electronic Engineers (IEEE, 2006) and Fellow of the International Federation for Automatic Control (IFAC, 2019). He was a co-recipient of the George S. Axelby Outstanding Paper Award for the Best Paper in IEEE Transactions on Automatic Control (2018). %He served as an Associate Editor for the journals Automatica, IEEE Transactions on Automatic Control, Systems and Control Letters, European Journal of Control and as a General Co-Chair of IEEE Conference on Decision and Control (CDC), Melbourne (2017). He currently serves as an Associate Editor for the IEEE Transactions on Control of Network Systems. \vspace{-0.5cm}
\end{IEEEbiography}

\begin{IEEEbiography}[{\includegraphics[width=1in,height=1.25in,clip,keepaspectratio]{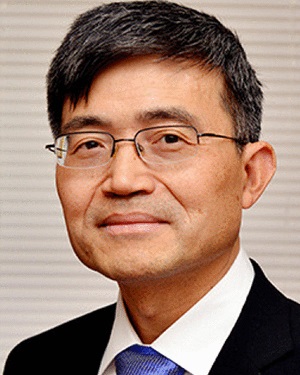}}]{Zhong-Ping Jiang}
(Fellow, IEEE) received the M.Sc. degree in statistics from the University of Paris XI, France, in 1989, and the Ph.D. degree in automatic control and mathematics from the Ecole des Mines de Paris (now, called ParisTech-Mines), France, in 1993, under the direction of Prof. Laurent Praly. Currently, he is a Professor of Electrical and Computer Engineering at the Tandon School of Engineering, New York University. His main research interests include stability theory, robust/adaptive/distributed nonlinear control, robust adaptive dynamic programming, reinforcement learning and their applications to information, mechanical and biological systems. He has served as Deputy Editor-in-Chief, Senior Editor and Associate Editor for numerous journals. Prof. Jiang is a Fellow of the IEEE, IFAC, and CAA, a foreign member of the Academia Europaea (Academy of Europe), and is among the Clarivate Analytics Highly Cited Researchers. In 2022, he received the Excellence in Research Award from the NYU Tandon School of Engineering.
\vspace{-0.5cm}
\end{IEEEbiography}

\end{document}